%Paper: Threshold $\theta \geq 2$ contact processes on homogeneous trees
%By: Luiz Renato Fontes and Roberto H. Schonmann

%Version submitted to AP Feb 2006 with minor corrections on Mar 2

\input amstex

\NoBlackBoxes
%\magnification=1200
\documentstyle{amsppt}
\document
\magnification=\magstep1
\hsize=5.3in
\vsize=7.4in
\TagsOnRight
\global\newcount\numsec\global\newcount\numfor
\global\newcount\numfig
\gdef\profonditastruttura{\dp\strutbox}
\def\senondefinito#1{\expandafter\ifx\csname#1\endcsname\relax}
\def\SIA #1,#2,#3 {\senondefinito{#1#2}
\expandafter\xdef\csname #1#2\endcsname{#3}\else
\write16{???? ma #1,#2 e' gia' stato definito !!!!} \fi}

\def\etichetta(#1){(\veroparagrafo.\veraformula)
\SIA e,#1,(\veroparagrafo.\veraformula)
 \global\advance\numfor by 1
% \write15{\string\FU (#1){\equ(#1)}}
 \write16{ EQ \equ(#1) == #1  }}

\def\letichetta(#1){\veroparagrafo.\veraformula
\SIA e,#1,{\veroparagrafo.\veraformula}
\global\advance\numfor by 1
% \write15{\string\FU (#1){\equ(#1)}}
 \write16{ Sta \equ(#1) == #1}}

\def\tetichetta(#1){\veroparagrafo.\veraformula %%%%copy four lines
\SIA e,#1,{(\veroparagrafo.\veraformula)}
\global\advance\numfor by 1
% \write15{\string\FU (#1){\equ(#1)}}
 \write16{ tag \equ(#1) == #1}}

\def \FU(#1)#2{\SIA fu,#1,#2 }

\def\etichettaa(#1){(A\veroparagrafo.\veraformula)
 \SIA e,#1,(A.\veroparagrafo.\veraformula)
\global\advance\numfor by 1
% \write15{\string\FU (#1){\equ(#1)}}
 \write16{ EQ \equ(#1) == #1  }}

\def\getichetta(#1){Fig. \verafigura
 \SIA e,#1,{\verafigura}
 \global\advance\numfig by 1
% \write15{\string\FU (#1){\equ(#1)}}
 \write16{ Fig. \equ(#1) ha simbolo  #1  }}

\newdimen\gwidth
\def\BOZZA{
\def\alato(##1){
 {\vtop to \profonditastruttura{\baselineskip
 \profonditastruttura\vss
 \rlap{\kern-\hsize\kern-1.2truecm{$\scriptstyle##1$}}}}}
\def\galato(##1){ \gwidth=\hsize \divide\gwidth by 2
 {\vtop to \profonditastruttura{\baselineskip
 \profonditastruttura\vss
 \rlap{\kern-\gwidth\kern-1.2truecm{$\scriptstyle##1$}}}}}
}
\def\alato(#1){}
\def\galato(#1){}
\def\veroparagrafo{\number\numsec}\def\veraformula{\number\numfor}
\def\verafigura{\number\numfig}
\def\Eq(#1){\eqno{\etichetta(#1)\alato(#1)}}
\def\teq(#1){\tag{\tetichetta(#1)\hskip-1.6truemm\alato(#1)}}  %%%%%this line
\def\eq(#1){\etichetta(#1)\alato(#1)}
\def\Eqa(#1){\eqno{\etichettaa(#1)\alato(#1)}}
\def\eqa(#1){\etichettaa(#1)\alato(#1)}
\def\eqv(#1){\senondefinito{fu#1}$\clubsuit$#1\else\csname fu#1\endcsname\fi}
\def\equ(#1){\senondefinito{e#1}\eqv(#1)\else\csname e#1\endcsname\fi}

%next six lines by paf (no responsibilities taken)
\def\Lemma(#1){Lemma \letichetta(#1)\hskip-1.6truemm}
\def\Theorem(#1){{Theorem \letichetta(#1)}\hskip-1.6truemm}
\def\Proposition(#1){{Proposition \letichetta(#1)}\hskip-1.6truemm}
\def\Corollary(#1){{Corollary \letichetta(#1)}\hskip-1.6truemm.}
\def\Remark(#1){{\noindent{\bf Remark \letichetta(#1)\hskip-1.6truemm.}}}

%%%%%%%%%%%%%%%%%% Numerazione verso il futuro ed eventuali paragrafi
%%%%%%%      precedenti non inseriti nel file da compilare
\def\include#1{
\openin13=#1.aux \ifeof13 \relax \else
\input #1.aux \closein13 \fi}
\openin14=\jobname.aux \ifeof14 \relax \else
\input \jobname.aux \closein14 \fi
%\openout15=\jobname.aux

%\BOZZA

\def\Bbbone{
     \kern.3em{\vrule height1.58ex width.6pt}%
     \kern.1em{\vrule height1.58ex width.6pt}%
     \kern-.392em\lower.12ex\hbox{\vrule height.4pt width.58em}%
     \kern-.4em\raise1.5ex\hbox{\vrule height.4pt width.34em}%
     \kern-.67em\raise.09ex\hbox{\' {}}%
     \kern-.56em\lower.13ex\hbox{\' {}}%
     \kern-.55em\lower.31ex\hbox{\' {}}%
     \kern.5em}

\def\({\left(}
\def\){\right)}

\def\P{\Bbb P}

\def\T{\Bbb T}

\def\Z{\Bbb Z}

\def\R{\Bbb R}

\def\one{1\kern-0.60em1}                %\Bbb 1}
\def\lone{1\kern-0.55em1}

\def\cqd{\hfill\vrule height .6em width .6em depth -.1em}

%\def\B{\Bbb B}
%%%%%%%%%%%%%%%%%%%%%%%%%%%%%%%%%%%%%%%%%%%%%%%%%%%%%%%%%%%%%%%%%%%
%\def\ssubset{\subset \subset}

%\def\G{\Cal{G}}

%\def\Ge{\Cal{G}_{\text{ext}}}

%\def\GA{\Cal{G}_{\text{A}}}

%\def\GAe{(\Cal{G}_{\text{A}})_{\text{ext}}}

%%%%%%%%%%%%%%%%%%%%%%%%%%%%%%%%%%%%%%%%%%%%%%%%%%%%%%%%%%%%%%%%%%%

%@@ Acknowledgments

\document

\topmatter
\title
Threshold $\theta \geq 2$ contact processes on homogeneous trees
\endtitle
\author Luiz Renato Fontes and Roberto H. Schonmann
   \footnote"*"{The work of L.R.F. was partially\
    supported by the Brazilian CNPq through grants 
    307978/2004-4 and 475833/2003-1, and by FAPESP through 
    grant 04/07276-2.
    The work of R.H.S. was partially supported by the American
    N.S.F. through grant DMS-0300672 .\newline}
%@@ support
\endauthor

\abstract
We study the threshold $\theta \geq 2$ contact process on a homogeneous tree
$\T_b$ of degree $\kappa = b + 1$, 
with infection parameter $\lambda \geq 0$ and 
started from a product measure with density $p$. 
The corresponding mean-field model displays a discontinuous transition at 
a critical point $\lambda_{\text{c}}^{\text{MF}}(\kappa,\theta)$ 
and for $\lambda \geq \lambda_{\text{c}}^{\text{MF}}(\kappa,\theta)$ it 
survives iff $p \geq p_{\text{c}}^{\text{MF}}(\kappa,\theta,\lambda)$,
where this critical density satisfies 
$0 < p_{\text{c}}^{\text{MF}}(\kappa,\theta,\lambda) < 1$, 
$ \lim_{\lambda \to \infty} 
p_{\text{c}}^{\text{MF}}(\kappa,\theta,\lambda) = 0$.
For large $b$, we show that the process on $\T_b$ has a qualitatively similar 
behavior when $\lambda$ is small, including the behavior at and close to 
the critical point $\lambda_{\text{c}}(\T_b,\theta)$. In contrast, for large
$\lambda$ the behavior of the process on $\T_b$ is qualitatively distinct
from that of the mean-field model in that the critical density has 
$p_{\text{c}}({\T}_b,\theta,\infty) :=
\lim_{\lambda \to \infty} p_{\text{c}}(\T_b,\theta,\lambda) > 0$.
We also show that $\lim_{b \to \infty} b \, \lambda_{\text{c}}(\T_b,\theta)
= \Phi_{\theta}$, where $1 < \Phi_2 < \Phi_3 < ...$, $\lim_{\theta \to \infty}
\Phi_{\theta} = \infty$, and 
$ 0 < \liminf_{b \to \infty} \  b^{\theta/(\theta-1)} \,
p_{\text{c}}({\T}_b,\theta,\infty)
\ \leq \
\limsup_{b \to \infty} \  b^{\theta/(\theta-1)} \,
p_{\text{c}}({\T}_b,\theta,\infty)
< \infty$.

\endabstract

\subjclass Primary 60K35
\endsubjclass
\keywords
Threshold Contact process, homogeneous trees, 
critical points, critical density, phase diagram, discontinuous transition
\endkeywords
%\affil IME-USP and Department of Mathematics, UCLA \endaffil
\affil USP and UCLA \endaffil
\endtopmatter
\leftheadtext{Luiz Renato Fontes and Roberto H. Schonmann}
\rightheadtext{Threshold contact processes}

%\subheading{1. INTRODUCTION AND RESULTS}
\subheading{1. Introduction and results}
\numsec=1
\numfor=1

\medskip

\noindent{\bf 1.1 Preliminaries.}
Threshold contact processes form a natural class of interacting particle
systems (see, e.g., [Lig1] for background to the area). 
They are most naturally defined on a graph or oriented graph $G = (V,E)$. 
In both cases $V$  is an arbitrary countable set, whose elements are called 
vertices or sites. When $G$ is a graph, 
$E \subset \{ \{v,u\} : v,u \in V \}$ is its set of edges, also called bonds.
When $G$ is an oriented graph, $E \subset V \times V$ is its set of 
oriented edges, also called oriented bonds.  
We denote the influence neighborhood of $v \in V$ in $G$ by 
$$
\Cal{N}_{G,v}  \  =  \ 
\left\{
\aligned 
& \{u \in V : \{v,u\} \in {E}\},
\quad \text{if $G$ is a graph},
\\
& \{u \in V : (v,u) \in {E}\},
\quad \text{if $G$ is an oriented graph}.
\endaligned
\right.
$$
The degree of the site $v \in V$ is the cardinality of $\Cal{N}_{G,v}$. 
The threshold $\theta$ contact process on ${G}$ is now defined as the 
Markov process on $\{0,1\}^V$ with flip rates at $v \in V$ at time 
$t \geq 0$ given by
\roster
\item"$\bullet$" 1 flips to 0 at rate 1.
\item"$\bullet$" 0 flips to 1 at rate $\lambda$ in case there are at least 
$\theta$ sites of $\Cal{N}_{G,v}$ in state 1 at time $t$, and at 
rate 0 otherwise. 
\endroster
The parameter $\lambda \geq 0$ is called the infection rate.
The state of the process at each site at each time is called the spin 
at that site at that time.
A spin 0 is interpreted as a vacant or healthy site, while a spin 1 
is interpreted as an occupied or infected site.
It is well known (see, e.g., Chapter I of [Lig1]) that such rates
define a unique Markov process. Note that the flip rates above are 
attractive (see Chapter III of [Lig1]), a property that has many 
consequences. 

When $\theta = 1$, the threshold contact process is easier to analyze,
among other reasons because it has an additive
dual process. (In this dual process 
infected sites become healthy at rate 1, and they infect 
simultaneously all their neighbors at rate $\lambda$.)
The behavior of the process is not expected then 
to be qualitatively different 
from that of the much studied (linear) contact process (see Chapter VI 
of [Lig1] and Part I of [Lig2]). 
For this reason we will focus in this
paper on the cases $\theta \geq 2$, which are more challenging and 
do present a different qualitative behavior. As we will explain in the 
next subsection, this difference in behavior is indicated by the associated
mean-field model.
In the end of this introduction we will review some results about threshold
$\theta \geq 2$ and related models, from [Toom], [DG], [BG], [Chen1], [Chen2].

For coupling purposes, 
it is convenient to construct the process using a system of Poisson marks. 
For this purpose, associate to each site in $V$ two independent Poisson
processes: one with rate 1, and one with rate $\lambda$. Mark the 
arrival times of the former with symbols $D$ (for ``down'') and those
of the latter with symbols $U$ (for ``up''). Make these Poisson 
processes independent from site to site. Use the marks now in the 
obvious way, to define the process: A spin 1 at site $v$
flips to 0 when it encounters a $D$ mark there; 
a spin 0 at site $v$ flips to 1 when it encounters an $U$ mark there
and at least $\theta$ neighbors of $v$ have spin $1$ at that time.
The probability space on which these Poisson processes are defined 
%(with probabilities denoted by $\P$) 
will be large enough to accommodate 
the process started from arbitrary initial configurations.

%Later on we will construct the process explicitly using flip 
%marks; the corresponding probability space (with probabilities
%denoted by $\P$) will be large enough to accommodate the process 
%started from arbitrary initial configurations.

We will denote by $(\eta^{\mu}_{{G},\theta,\lambda;t})_{t \geq 0}$ 
the process started from a random distribution picked according 
to law $\mu$ at time 0. When $\mu$ is product measure with density 
$p$ we will use the notation $(\eta^{p}_{{G},\theta,\lambda;t})$. When 
there is no risk of confusion, 
${G}$, $\theta$ and $\lambda$ may be omitted from the notation.

The point mass on the configuration with all sites in state $i \in \{0,1\}$
will be denoted $\delta_{G,i}$. The distribution 
$\delta_{G,0}$ is trivially invariant for the threshold $\theta$ 
contact process, when $\theta > 0$.
By attractivity, $\eta^{1}_{{G},\theta,\lambda;t} \Rightarrow 
\nu_{{G},\theta,\lambda}$, as $t \to \infty$,  
where $\Rightarrow$ denotes convergence in distribution, and 
$\nu_{{G},\theta,\lambda}$ is called the upper invariant measure. 

We say that the process started from the distribution $\mu$ dies out
when $\eta^{\mu}_{{G},\theta,\lambda;t} \Rightarrow \delta_{G,0}$, 
as $t \to \infty$. When this happens for every $\mu$, we simply say 
that the process dies out. Attractivity implies that the process dies
out precisely when $\nu_{{G},\theta,\lambda} = \delta_{G,0}$.
When the process does not die out, we will say that it survives.

For $v \in V$, set 
$$
\rho^{\mu}_{{G},\theta,\lambda;t}(v)
\ = \
\P(\eta^{\mu}_{{G},\theta,\lambda;t}(v) = 1).
$$
We will use for $\rho^{\mu}_{G,\theta,\lambda;t}(v)$ the same 
conventions on notation as for $\eta^{\mu}_{G,\theta,\lambda;t}$.
Also, when $G$ is such that 
$\rho^{\mu}_{G,\theta,\lambda;t}(v)$ does not depend on $v$, we will omit 
$v$ from this notation. 
The critical point for the threshold $\theta$ contact process on 
${G}$ is defined by
$$
\lambda_{\text{c}}({G},\theta) 
\ = \ 
\sup\{\lambda : \nu_{{G},\theta,\lambda} = \delta_{G,0}\} 
\ = \
\sup\{\lambda : \text{for each $v \in V$,
$\rho^{1}_{{G},\theta,\lambda;t}(v) \to 0$ as $t \to \infty$}\}.
$$
Clearly, the convergence of $\rho^{1}_{{G},\theta,\lambda;t}(v)$
to 0 cannot be faster then exponential. Explicitly:
$$
\rho^{1}_{{G},\theta,\lambda;t}(v) 
\ \geq  \
\P (\text{there is no $D$ mark at $v$ from time 0 to $t$})
\ \geq  \
e^{-t}, 
%e^{-t}, \quad \text{ for $t \geq 0$}.
\Eq(explowerbd)
$$
for $t \geq 0$.
It is natural to define
$$
\lambda_{\text{exp}}({G},\theta) 
\ = \
\sup\{\lambda : \text{for each $v \in V$,
$\rho^{1}_{{G},\theta,\lambda;t}(v) \to 0$ 
exponentially fast as $t \to \infty$}\}.
$$
Obviously 
$$
\lambda_{\text{exp}}({G},\theta)
\ \leq \ 
\lambda_{\text{c}}({G},\theta),
$$
and it is interesting to decide when equality holds.

Even when the process survives, it may happen that for small $p>0$, 
$\eta^{p}_{{G},\theta,\lambda;t} \Rightarrow
\delta_{G,0}$, as $t \to \infty$. 
%(see [DG], [Che1], [Che2], [FSS]). 
By attractivity, if this happens for some value of $p$, it will
also happen for smaller values of $p$. This leads to the following
definitions:
$$
\align
p_{\text{c}}({G},\theta,\lambda)
& \ = \
\sup\{p \in [0,1] : \eta^{p}_{{G},\theta,\lambda;t} \Rightarrow
\delta_{G,0}, \ \text{as} \ t \to \infty \}  \\
& \ = \
\sup\{p \in [0,1]: \text{for each $v \in V$,
$\rho^{p}_{{G},\theta,\lambda;t}(v) \to 0$ as $t \to \infty$}\},
\endalign
$$
$$
p_{\text{exp}}({G},\theta,\lambda)
\ = \
\sup\{p \in [0,1]: \text{for each $v \in V$,
$\rho^{p}_{{G},\theta,\lambda;t}(v) \to 0$ 
exponentially fast as $t \to \infty$}\}.
$$
As above, obviously
$$
p_{\text{exp}}({G},\theta,\lambda)
\ \leq \
p_{\text{c}}({G},\theta,\lambda),
$$
and again it is interesting to decide when equality holds.

When $(\eta^{p}_{{G},\theta,\lambda;t})$ survives, it is natural to 
ask if it converges in distribution as $t \to \infty$. This is in 
general a difficult question, but the following provides 
a partial answer and, as a by-product, an estimate on 
$p_{\text{c}}({G},\theta,\lambda)$. 

\proclaim{Proposition 1} For any $G$, $\theta$ and $\lambda$, 
for $p \geq \lambda/(\lambda + 1)$,
$$
\eta^{p}_{G,\theta,\lambda;t}  \Rightarrow
\nu_{G,\theta,\lambda},
\quad \text{as $t \to \infty$}.
\Eq(Proposition11)
$$
Therefore, if the process survives for a certain value of $\lambda$,
then
$$
p_{\text{c}}(G,\theta,\lambda)
\leq 
\frac{\lambda}{\lambda + 1}.
\Eq(Proposition12)
$$
\endproclaim
\demo{Proof}
Let $\beta_{G,p}$ be the product measure with density $p$.
For any $\theta \geq 0$, 
the threshold $\theta$ contact process is stochastically
dominated by the threshold 0 contact process. This latter process is simply
an independent flip process whose unique invariant distribution is
$\beta_{G,\lambda / (\lambda + 1)}$.  
Hence, for $p \geq \lambda / (\lambda + 1)$,
$$
\nu_{G,\theta,\lambda} \ \leq \ \beta_{G,\lambda / (\lambda + 1)} \ \leq 
\ \beta_{G,p} \ \leq \ \delta_{G,1}, \quad \text{stochastically}.
$$
But as $t \to \infty$, we know that 
$\eta_{G,\theta,\lambda;t}^{\nu_{G,\theta,\lambda}} 
\Rightarrow \nu_{G,\theta,\lambda}$
and $\eta_{G,\theta,\lambda;t}^1 \Rightarrow \nu_{G,\theta,\lambda}$. 
Therefore \equ(Proposition11) also holds.
%$$ \eta_t^p \Rightarrow \nu.  $$
\cqd
\enddemo
\noindent
{\it Remark on case $\theta = 1$:} In that case, duality can easily be
used to show that in Proposition 1, \equ(Proposition11) holds for every 
$p > 0$ and hence \equ(Proposition12) can be replaced by 
$p_{\text{c}}(G,\theta,\lambda) = 0$.

\medskip

\noindent{\bf 1.2 The mean-field model.}
When $G$ is regular of degree $\kappa$ (i.e., each site has 
$\kappa$ neighbors), 
its is natural to compare the evolution of the threshold contact process on
$G$ with a corresponding ``mean-field'' evolution. 
By this we mean the evolution of a deterministic density 
$(\rho^{\text{MF},p}_{\kappa,\theta,\lambda;t})_{t \geq 0}$,
which is governed by
$$
\frac{d}{dt} \rho^{\text{MF},p}_{\kappa,\theta,\lambda;t}
\ = \
-\rho^{\text{MF},p}_{\kappa,\theta,\lambda;t}
\ + \
\lambda \, 
\left(1-\rho^{\text{MF},p}_{\kappa,\theta,\lambda;t}\right) 
\,
\text{Bin}\left(
\kappa,\rho^{\text{MF},p}_{\kappa,\theta,\lambda;t},\theta \right),
\Eq(MF1)
$$
with $\rho^{\text{MF},p}_{\kappa,\theta,\lambda;0} = p$,
and where
$$
\text{Bin}(\kappa,x,\theta)
\ = \
\sum_{i=\theta}^{\kappa} \ \left(
\aligned  & \kappa \\ & i \endaligned
\right) x^i (1-x)^{\kappa-i} 
$$
is the probability that a binomial
random variable with $\kappa$ attempts and probability of success $x$
is larger than or equal to $\theta$.

To explain the origin of this mean-field evolution and its relationship 
with the threshold contact process on $G$, we observe that from 
the definition of that process, for each $v \in V$, 
$$
\frac{d}{dt} \rho^{p}_{G,\theta,\lambda;t} (v)
\ = \
- \rho^{p}_{G,\theta,\lambda;t} (v)
\
+ \lambda \,
\, \P\left(\eta^{p}_{G,\theta,\lambda;t}(v) = 0,
%\ \left|\left\{
\ \#\left\{
u \in \Cal{N}_{G,v} : \eta^{p}_{G,\theta,\lambda;t}(u) = 1
%\right\} \right| \geq \theta\right).
\right\}  \geq \theta\right).
%\Eq(generator)
$$
One then obtains \equ(MF1) if one pretends that the $\kappa + 1$ random
variables $\eta^{p}_{G,\theta,\lambda;t}(v)$, 
$\eta^{p}_{G,\theta,\lambda;t}(u)$, $u \in \Cal{N}_{G,v}$, 
are i.i.d., 
with common density $\rho^{\text{MF},p}_{\kappa,\theta,\lambda;t}$.

Trivially, $\rho^{\text{MF},0}_{\kappa,\theta,\lambda;t} = 0$, 
for all $t \geq 0$.
When $p \in (0,1]$, it clearly follows from \equ(MF1) that 
$\rho^{\text{MF},p}_{\kappa,\theta,\lambda;t} \geq p e^{-t} > 0$, 
for all $t \geq 0$. 
It will then be convenient to rewrite \equ(MF1) as
$$
\frac{d}{dt} \log (\rho^{\text{MF},p}_{\kappa,\theta,\lambda;t})
=
H(\kappa,\theta, \lambda;\rho^{\text{MF},p}_{\kappa,\theta,\lambda;t}),
\Eq(MF2)
$$
where
$$
H(\kappa,\theta, \lambda;x)
=
-1
+
\lambda \frac{1-x}{x} \text{Bin}(\kappa,x,\theta),
\qquad x \in (0,1].
$$
Note that when $\theta \geq 2$, then
$0 \leq \text{Bin}(\kappa,x,\theta) \leq (\kappa^2/2)x^2$, and hence
$$
\lim_{x \searrow 0} H(\kappa,\theta, \lambda;x)
=
-1.
\Eq(xto0H)
$$
Also important are the elementary facts that $H(\kappa,\theta, \lambda;x)$
is continuous in $x$ and in $\lambda$, 
$H(\kappa,\theta, \lambda;1) = -1$, 
$H(\kappa,\theta,0;x) = -1$, 
and, provided $\kappa \geq \theta$,  
for each $x \in (0,1)$, $H(\kappa,\theta, \lambda;x)$
is strictly increasing in $\lambda$, with 
$\lim_{\lambda \to \infty} H(\kappa,\theta, \lambda;x) = \infty$.

%From these facts it follows that there is a critical point 
%$\lambda_{\text{c}}^{\text{MF}}(\kappa,\theta)$ separating a low $\lambda$
%regime in which $\sup_{x\in(0,1]}F(\kappa,\theta, \lambda;x) < 0$,
%from a large $\lambda$ regime in which 

These facts motivate the definition of the critical point
$$
\lambda_{\text{c}}^{\text{MF}}(\kappa,\theta)
=
\sup\{\lambda \geq 0: \sup_{x\in(0,1]}H(\kappa,\theta, \lambda;x) < 0 \},
$$
and the critical density
$$
p_{\text{c}}^{\text{MF}}(\kappa,\theta,\lambda)
=
\inf \{x \in (0,1] : H(\kappa,\theta, \lambda;x) \geq 0 \}.
%\sup \{x \in (0,1] : H(\kappa,\theta, \lambda;y) < 0
%\text{ \ for all $y \in (0,x]$} \}.
$$
Note that the facts above imply that, when $\kappa \geq \theta \geq 2$,  
$$
0 < \lambda_{\text{c}}^{\text{MF}}(\kappa,\theta) < \infty, 
$$
and 
$$
0 < p_{\text{c}}^{\text{MF}}(\kappa,\theta,\lambda) < 1 \quad \text{for} 
\ \lambda \geq \lambda_{\text{c}}^{\text{MF}}(\kappa,\theta).
$$
Moreover 
$p_{\text{c}}^{\text{MF}}(\kappa,\theta,\lambda)$ is strictly 
decreasing in $\lambda \geq \lambda_{\text{c}}^{\text{MF}}(\kappa,\theta)$,
with 
$$
\lim_{\lambda \to \infty} p_{\text{c}}^{\text{MF}}(\kappa,\theta,\lambda)
= 0.
\Eq(pcinftyMF)
$$
Define 
$$
\Cal{D}^{\text{MF}}(\kappa,\theta) 
\ = \
\left\{ 
(\lambda,p) \in [0,\infty] \times [0,1] : 
\lambda < \lambda_{\text{c}}^{\text{MF}}(\kappa,\theta)
\text{ or }
p < p_{\text{c}}^{\text{MF}}(\kappa,\theta,\lambda)
\right\}.
$$
\proclaim{Proposition 2}
For every $\kappa \geq \theta \geq 2$, the following dichotomy holds.

In case $(\lambda,p) \in (\Cal{D}^{\text{MF}}(\kappa,\theta))^c$, 
$$
\rho^{\text{MF},p}_{\kappa, \theta,\lambda;t}
\ \geq \
p_{\text{c}}^{\text{MF}}(\kappa,\theta,\lambda)
\ > \ 0,
\text{ \ for all $t \geq 0$}.
\Eq(Proposition2a)
$$
In case $(\lambda,p) \in \Cal{D}^{\text{MF}}(\kappa,\theta)$,
then for some $C \in (0,\infty)$,
$$
p \, e^{-t} 
\leq \
\rho^{\text{MF},p}_{\kappa, \theta,\lambda;t}
\ \leq \
C \, e^{-t}, \text{ \ for all $t \geq 0$}.
\Eq(Proposition2b)
$$
\endproclaim
In particular, a discontinuous transition happens at
$\lambda_{\text{c}}^{\text{MF}}(\kappa,\theta)$.
\demo{Proof}
The only statement that requires explanation is the upper bound in 
\equ(Proposition2b). To prove it, note first that \equ(MF2) implies
that when $(\lambda,p) \in \Cal{D}^{\text{MF}}(\kappa,\theta)$
%$p < p_{\text{c}}^{\text{MF}}(\kappa,\theta,\lambda)$
we have $\lim_{t \to \infty} \rho^{\text{MF},p}_{\kappa, \theta,\lambda;t}
= 0$.
Using \equ(MF2) a second time now, this time in combination with 
\equ(xto0H), shows that for any $\epsilon > 0$, there is $C' \in (0,\infty)$
such that 
$$
\rho^{\text{MF},p}_{\kappa, \theta,\lambda;t}
\ \leq \
C' \, e^{-(1-\epsilon)t}, \text{ \ for all $t \geq 0$}.
$$
Using this estimate in combination with \equ(MF1) now yields, since 
$\theta \geq 2$, 
$$
\frac{d}{dt} \rho^{\text{MF},p}_{\kappa,\theta,\lambda;t}
\ \leq  \
-\rho^{\text{MF},p}_{\kappa,\theta,\lambda;t}
\ + \
\lambda \, 
\frac{\kappa^2}{2} \, (C')^2 \, e^{-2(1-\epsilon)t}.
$$
%Multiplying by $e^t$ and integrating,
Multiplying by $e^t$, we obtain
$$
\frac{d}{dt} \left( 
e^t \, \rho^{\text{MF},p}_{\kappa,\theta,\lambda;t}
\right)
\ \leq \ 
\lambda \, \frac{\kappa^2}{2} \, (C')^2 \, e^{-(1-2\epsilon)t}.
$$
Supposing $\epsilon < 1/2$, integration in $t$ from 0 to $s$ yields
$$
e^s \, \rho^{\text{MF},p}_{\kappa,\theta,\lambda;s}
%\ - \ \rho^{\text{MF},p}_{\kappa,\theta,\lambda;0}
\ -  \ p
\ \leq \
\int_0^{s}
\lambda \, \frac{\kappa^2}{2} \, (C')^2 \, e^{-(1-2\epsilon)t} \, dt
\ \leq \
\int_0^{\infty} 
\lambda \, \frac{\kappa^2}{2} \, (C')^2 \, e^{-(1-2\epsilon)t} \, dt
\ = \ C''\ < \ \infty.
$$
\equ(Proposition2b) follows, with $C = C'' + p$. 
\cqd
\enddemo

It is easy to see that $\lambda_{\text{c}}^{\text{MF}}(\kappa,\theta)
\to 0$, as $\kappa \to \infty$. The rate at which this convergence 
occurs is also easily identified from standard facts about convergence 
of binomial distributions to Poisson distributions.
For this purpose extend the definition of 
$\text{Bin}(\kappa,x,\theta)$ to be 1 when $x > 1$, and observe that 
straightforward computations then yield
$$
\text{Bin}(\kappa,\gamma/\kappa,\theta) 
\ \to \ \text{Poisson}(\gamma,\theta) 
\ := \ \sum_{i\geq \theta} e^{-\gamma} \frac{\gamma^i}{i!},
\quad \text{as $\kappa \to \infty$},
\Eq(bintopoisson)
$$
uniformly in $\gamma > 0$. 
The corresponding limit for the function $H(\kappa,\theta,\lambda;x)$ is
$$
H(\kappa,\theta,\phi/\kappa;\gamma/\kappa) \ \to \ 
-1 \ + \ \phi \ \frac{\text{Poisson}(\gamma,\theta)}{\gamma}, 
\quad \text{as $\kappa \to \infty$},
$$
uniformly in $\gamma > 0$, for each $\phi > 0$.
From this it is easy to derive 
$$
\lim_{\kappa \to \infty} \,
\kappa \, \lambda_{\text{c}}^{\text{MF}}(\kappa,\theta) 
\ = \ 
\Phi_{\theta} 
\ := \
\ \inf_{\gamma > 0} \
\frac{\gamma}{\text{Poisson}(\gamma,\theta)}.
\Eq(kappalambdac)
$$
The constants $\Phi_{\theta}$ can easily be shown to satisfy
$$
1 < \Phi_2 < \Phi_3 < ...
\qquad \text{and} \qquad
\lim_{\theta \to \infty} \Phi_{\theta} = \infty.
$$
{\it Remark on case $\theta = 1$:} In this case, in contrast to 
\equ(xto0H), we have 
$\lim_{x \searrow 0} H(\kappa,\theta, \lambda;x) = -1 + \lambda \kappa$.
It is an instructive exercise to use this fact and the bound 
$H(\kappa,\theta, \lambda;x) < -1 + \lambda \kappa$, for $x > 0$, 
to analyze the behavior of the mean-field model in this case. In 
contrast to Proposition 2, one finds a continuous transition at 
$\lambda_{\text{c}}^{\text{MF}}(\kappa,1) = 1/\kappa$, with 
$p_{\text{c}}^{\text{MF}}(\kappa,1,\lambda) = 0$ for all 
$\lambda > \lambda_{\text{c}}^{\text{MF}}(\kappa,1)$.
The analogue of \equ(kappalambdac) is also true, with $\Phi_1 = 1$.

\medskip

\noindent{\bf 1.3 Results for the process on homogeneous trees.}
We will study threshold contact processes on the homogeneous tree
$\T_b$, of degree $b + 1$. 

Our first result provides conditions for survival of the process
on $\T_b$ based on the survival of the mean-field model with 
$\kappa = b$ (note: not $\kappa = b+1$) and
$\lambda/(\lambda + 1)$ in place of $\lambda$.  

\proclaim{Theorem 1} For arbitrary $b \geq 2$ and $\theta \geq 2$,
if $\lambda_{\text{c}}^{\text{MF}}(b,\theta) < 1$, then
$$
\lambda_{\text{c}} ({\T}_b,\theta)
\ \leq \ \frac{\lambda_{\text{c}}^{\text{MF}}(b,\theta)}{
1 - \lambda_{\text{c}}^{\text{MF}}(b,\theta)} 
\ < \ \infty,
\Eq(Thm11)
$$
and for arbitrary $\lambda > \lambda_{\text{c}}^{\text{MF}}(b,\theta)/
(1 - \lambda_{\text{c}}^{\text{MF}}(b,\theta))$,
$$
p_{\text{c}}({\T}_b,\theta,\lambda)
\ \leq  \ p_{\text{c}}^{\text{MF}}(b,\theta,\lambda/(\lambda + 1)) \ < \ 1.
\Eq(Thm12)
$$
\endproclaim
Since we know that 
$\lim_{b \to \infty} \lambda_{\text{c}}^{\text{MF}}(b,\theta) = 0$, 
Theorem 1 contains meaningful statements when $b$ is large. Moreover, 
combining this theorem with \equ(kappalambdac), we learn that 
$$
\limsup_{b \to \infty} 
\ b \,\lambda_{\text{c}} ({\T}_b,\theta)
\ \leq \ \Phi_{\theta}.
\Eq(Thm31<)
$$
This result will be sharpened in Theorem 3 below.
In our approach to the proof of Theorem 3, we will prove first the somewhat 
technical Theorem 2 below. 
Note that thanks to \equ(Thm31<), this theorem covers the behavior near
$\lambda_{\text{c}} ({\T}_b,\theta)$, when $b$ is large.
This theorem should be compared to Proposition 2. 
\proclaim{Theorem 2} For arbitrary $\theta \geq 2$ and $A \in (0,\infty)$,
there are $b_0,\delta \in (0,\infty)$
such that if $b \geq b_0$ and $\lambda \leq A/b$, then for every $p \in [0,1]$
the following dichotomy holds. Either
$$
\liminf_{T \to \infty} \, \frac{1}{T}
\int_0^T \, \rho^{p}_{\T_b,\theta,\lambda;t} \, dt
\ \geq \
\frac{\delta}{b},
\Eq(Thm41)
$$
or, for some $C \in (0,\infty)$,
$$
\rho^{p}_{\T_b,\theta,\lambda;t}
\ \leq \
C \, e^{-t}, \text{ \ for all $t \geq 0$}.
\Eq(Thm42)
$$
Moreover, for $b \geq b_0$ the set
$$
\Cal{D}_A(\T_b,\theta) 
\ = \ \{(\lambda,p) \in [0,A/b] \times [0,1] :
\text{alternative \equ(Thm42) holds}  \}
$$
is an open subset of $[0,A/b] \times [0,1]$ in the relative topology
induced by the Euclidean topology of $\R^2$. 
\endproclaim
Note that alternative \equ(Thm41) implies that the process survives when 
it starts with density $p$. 
For $p=1$, \equ(Thm41) is equivalent to each one of the statements
$$
\rho^{1}_{\T_b,\theta,\lambda;t} \ \geq \ \frac{\delta}{b},
\text{ \ for all $t \geq 0$},
$$
and
$$
\rho^{1}_{\T_b,\theta,\lambda;\infty} 
\ := \ \lim_{t \to \infty} \rho^{1}_{\T_b,\theta,\lambda;t}
\ \geq \ \frac{\delta}{b}.
$$

The following is immediate from Theorem 2, \equ(Thm31<) and 
\equ(Proposition12).
\proclaim{Corollary 1}
For each $\theta \geq 2$, the following statements hold when $b$ is large.
$$
0 \ < \
\lambda_{\text{exp}} ({\T}_b,\theta) 
\ = \
\lambda_{\text{c}} ({\T}_b,\theta)
\ < \ \infty,
$$
and the process survives at this critical point.
Moreover for $\lambda \geq \lambda_{\text{c}} ({\T}_b,\theta)$ close to 
this critical point, 
%$\lambda_{\text{c}} ({\T}_b,\theta)$,
%Moreover for $\lambda$ close to and not smaller than
%$\lambda_{\text{c}} ({\T}_b,\theta)$,
$$
0 \ < \
p_{\text{exp}}({\T}_b,\theta,\lambda) 
\ = \ 
p_{\text{c}}({\T}_b,\theta,\lambda)
\ < \ 1,
$$
and the process started from this critical density
$p_{\text{c}}({\T}_b,\theta,\lambda)$ survives.
\endproclaim
%Combining Corollary 1 and Proposition 1, we see
%that the process on $\T_b$, with large $b$,
%shares one more qualitative feature with the mean-field model:
Note that in particular, under the conditions in Corollary 1, 
$$
p_{\text{c}}({\T}_b,\theta,\lambda_{\text{c}} ({\T}_b,\theta))
\ < \ 1.
$$

%%%%%%%%%%%%%%%%%%%%%%%%%%%%%%%%%%%%%%%%%%%%%%%%%%%%%%%

Theorem 2 and Corollary 1 show qualitative similarities
between the behavior of the threshold contact process on $\T_b$ and the
corresponding mean-field model, when b is large and $\lambda$ is small.
The next theorem shows a related quantitative similarity.
\proclaim{Theorem 3} For arbitrary $\theta \geq 2$,
$$
\lim_{b \to \infty} b \, \lambda_{\text{exp}} ({\T}_b,\theta)
\ = \
\lim_{b \to \infty} b \, \lambda_{\text{c}} ({\T}_b,\theta)
\ = \
\Phi_{\theta}.
\Eq(Thm31)
$$
\endproclaim
\noindent
{\it Remark on case $\theta = 1$:} In that case, duality can easily be
used to show that \equ(Thm31) also holds, with $\Phi_1 = 1$.

In contrast to the results above,
for large values of $\lambda$, the process on $\T_b$ and 
the mean-field model behave differently, as the comparison between
the following theorem and \equ(pcinftyMF) shows. 
Set
$$
p_{\text{exp}}({\T}_b,\theta,\infty)
\ = \
\lim_{\lambda \to \infty} p_{\text{exp}}({\T}_b,\theta,\lambda),
\qquad
p_{\text{c}}({\T}_b,\theta,\infty)
\ = \
\lim_{\lambda \to \infty} p_{\text{c}}({\T}_b,\theta,\lambda).
$$
%$$
%p_{\text{c}}({\T}_b,\theta,\infty)
%\ = \
%\lim_{\lambda \to \infty} p_{\text{c}}({\T}_b,\theta,\lambda).
%$$

\proclaim{Theorem 4} For arbitrary $b \geq 2$ and $\theta \geq 2$,
$$
p_{\text{c}}({\T}_b,\theta,\infty)
\ \geq \
p_{\text{exp}}({\T}_b,\theta,\infty)
\ > \ 0.
$$
\endproclaim

We do not know if the critical densities 
$p_{\text{c}}({\T}_b,\theta,\infty)$ 
and $p_{\text{exp}}({\T}_b,\theta,\infty)$ are identical to each other, 
but the next theorem shows that at least they display similar behavior as
$b \to \infty$.
\proclaim{Theorem 5} For arbitrary $\theta \geq 2$,
$$
0 <
\liminf_{b \to \infty} \  b^{\theta/(\theta-1)} \,
p_{\text{exp}}({\T}_b,\theta,\infty)
\ \leq \
\limsup_{b \to \infty} \  b^{\theta/(\theta-1)} \,
p_{\text{c}}({\T}_b,\theta,\infty)
< \infty.
\Eq(Thm32)
$$
\endproclaim

\medskip

\noindent{\bf 1.4 Results for the process on oriented homogeneous trees.}
As a tool, for comparison purposes and for its own sake, we will also
study the threshold contact process on an oriented graph $\vec\T_b$, 
obtained from $\T_b$ in a fashion described next. 

First we introduce some notation. 
We embed a copy of $\Z$ in $\T_b$ and use the notation 
$\Cal{L} = \Cal{L}_0$ to 
denote the set of vertices of $\T_b$ covered by this embedding. 
Abusing notation, we will denote the elements of $\Cal{L}$ by the 
names of the elements of $\Z$ that they represent in this embedding. 
We will also refer to site 0 as the root of $\T_b$.
Define, inductively in $n \geq 1$, $\Cal{L}_n$ as the set of vertices 
which are neighbors to some vertex in $\Cal{L}_{n-1}$
and are not in $\cup_{i = 0}^{n-1} \Cal{L}_{i}$. 
%and are not in $\Cal{L}_{n-2}$. 

With this notation, 
include $(v,u)$ in the set of oriented edges of $\vec\T_b$ if
$v,u \in \Cal{L}$ and $u = v+1$, or if 
$v \in \Cal{L}_{n-1}$ and $u \in \Cal{L}_n$, for some $n \geq 1$.
%the path from $u$ to $\Cal{L}$
%in $\T_b$ contains the edge $\{v,u\}$. Note that in $T_b$ each
%site $v$ has degree $|\Cal{N}_v| = b+1$, while in $\vec{T}_b$ each
%site $v$ has $|\Cal{N}_v| = b$.

Theorems 1 to 5 have analogues for the threshold $\theta$ contact 
process on $\vec\T_b$. 
Theorem 1 for $\T_b$ is actually a corollary to the same statement for 
$\vec{\T}_b$. Theorem 2 and Corollary 1 admit much  
stronger versions for $\vec{\T}_b$; those are stated as Theorem 6 and 
Corollary 2 below.
The analogues of Theorems 3, 4 and 5 for $\vec{\T}_b$ are also true
and can either be obtained
from the corresponding results for $\T_b$, or can more easily be proved 
directly.

\proclaim{Theorem 6} For arbitrary $\theta \geq 2$ and $b \geq 2$,
the following dichotomy holds
for every $\lambda \geq 0$ and $p \in [0,1]$.
Either
$$
\rho^{p}_{\vec{\T}_b,\theta,\lambda;t}
\ \geq \
p_{\text{c}}^{\text{MF}}(b,\theta,\lambda),
\text{ \ for all $t \geq 0$},
\Eq(Thm41a)
$$
or, for some $C \in (0,\infty)$,
$$
\rho^{p}_{\vec{\T}_b,\theta,\lambda;t}
\ \leq \
C \, e^{-t}, \text{ \ for all $t \geq 0$}.
\Eq(Thm42a)
$$
Moreover, for every $\theta \geq 2$ and $b \geq 2$,
the set
$$
\Cal{D}(\vec{\T}_b,\theta) 
\ = \ \{(\lambda,p) \in [0,\infty) \times [0,1] :
\text{alternative \equ(Thm42a) holds}  \}
$$
is an open subset of $[0,\infty) \times [0,1]$ in the relative topology
induced by the Euclidean topology of $\R^2$.
\endproclaim

\proclaim{Corollary 2}
For each $\theta \geq 2$
and $b \geq 2$ for which $\lambda_{\text{exp}} (\vec{\T}_b,\theta) < \infty$, 
the following statements hold.
$$
\lambda_{\text{c}}^{\text{MF}}(b,\theta) 
\ \leq \
\lambda_{\text{exp}} (\vec{\T}_b,\theta) 
\ = \ 
\lambda_{\text{c}} (\vec{\T}_b,\theta),
$$
and the process survives at this critical point.
Moreover for $\lambda \geq \lambda_{\text{c}} (\vec{\T}_b,\theta)$,
$$
p_{\text{c}}^{\text{MF}}(b,\theta,\lambda)
\ \leq \
p_{\text{exp}}(\vec{\T}_b,\theta,\lambda) 
\ = \ 
p_{\text{c}}(\vec{\T}_b,\theta,\lambda) 
\ < \ 
1,
$$
and the process started from this critical density
$p_{\text{c}}(\vec{\T}_b,\theta,\lambda)$ survives.
\endproclaim

We do not know if the first inequality in each display in Corollary 2
holds for the threshold contact process on $\T_b$. But it is known that 
the second of these does not hold for the process on 
$\Z^d$, $d \geq 3$, with $\theta = 2$, as reviewed in the next subsection.

\medskip

\noindent{\bf 1.5 Related previous results.} 
Threshold contact processes with $\theta = 2$ and closely related models
have been studied in [Toom], [DG], [BG], [Chen1], [Chen2], sometimes under
the name ``sexual contact process'', and mostly on $\Z^d$. 
In [Toom] discrete time versions were studied, 
and contour arguments were used to show survival. In [DG] these
contour methods were adapted to continuous time; their main result can be 
stated as follows using our terminology. Let $\vec{\Z}^2$ be the oriented 
graph obtained from $\Z^2$ by setting 
$\Cal{N}_{\vec{\Z}^2,v} = \{v+(1,0), \, v + (0,1) \}$. It is 
proved in [DG] that $\lambda_{\text{c}}(\vec{\Z}^2,2) < \infty$. 
Note that, by stochastic domination, this implies that 
$\lambda_{\text{c}}({\Z}^d,2) < \infty$, for $d \geq 2$.
In [BG] a renormalization procedure was introduced, which can replace 
the contour methods in proving survival. In [Chen1] and [Chen2] 
continuous time models which can be seen as modified 
threshold $\theta = 2$ contact processes on $\Z^d$ were studied. 
In one of these modified models
the flip rates at $v \in \Z^d$ at time $t \geq 0$ are given by
\roster
\item"$\bullet$" 1 flips to 0 at rate 1.
\item"$\bullet$" 0 flips to 1 at rate $\lambda$ in case there are at least
$2$ sites of $\Cal{N}_{\Z^d,v}$ that are separated from each other
by Euclidean distance $\sqrt{2}$ and are in state 1 at time $t$.
%and at rate 0 otherwise.
\endroster
The most important result from [Chen1], [Chen2] in connection to the 
current paper is the fact that for this modified model in $d \geq 3$, when 
$\lambda$ is large, survival occurs starting from any positive density $p$.
Since the threshold 2 contact process dominates that modified 
model, we learn that when $d \geq 3$ and $\lambda$ is large,
$p_{\text{c}}(\Z^d, 2, \lambda) = 0$. 
This means that the qualitative behavior of the $\theta = 2$
threshold contact process on $\Z^d$, $d \geq 3$, deviates from that of
the corresponding mean-field models (for which 
$p_{\text{c}}^{\text{MF}}(2d, 2, \lambda) > 0$, for all $\lambda> 0$), 
but this deviation is in the ``opposite direction'' of the deviation 
observed in the corresponding models on homogeneous trees
(for which, contrary to the 
mean-field model, $\lim_{\lambda \to \infty}
p_{\text{c}}(\T_{2d-1}, 2, \lambda) > 0$).

The only results that we are aware of for threshold $\theta \geq 2$ 
contact processes on trees are the following ones, from [DG]. There the
authors consider the model with $\theta = 2$ on $\vec{\T}_2$. They 
state that the contour methods used in that paper can be used to prove
that this process survives when $\lambda$ is large. They then prove, 
using this result, that the transition at $\lambda_{\text{c}}(\vec{\T}_2,2)$
is discontinuous (our Corollary 2 extends this result).
%Our Theorem 6 and Corollary 3 extend this result. 
%One of our motivations in the current paper is also to extend this 
%result to $\T_b$. 

\medskip

\noindent{\bf 1.6 Organization of the paper.} In Section 2 we prove the results
about the threshold contact process on $\vec{\T}_b$ (obtaining Theorem 1
as a corollary). In Section 3 we prove the results about 
the threshold contact process on ${\T}_b$ when $b$ is large and $\lambda$
is small, namely, Theorems 2 and 3. In Section 4 we prove 
the results about the threshold contact process on ${\T}_b$
when $\lambda$ is large, namely, Theorems 4 and 5; for this purpose
the bootstrap percolation model will be introduced as a tool. 

\

\subheading{2. Comparison between the model on $\vec{\T}_b$ and the 
mean-field model}
\numsec=2
\numfor=1
In this section we will prove the analogue of Theorem 1 for $\vec{\T}_b$
and Theorem 6. Both theorems result from a fairly direct comparison with the 
mean-field model. This comparison is based on writing down, for 
an arbitrary site $v$, the differential equation
$$
\frac{d}{dt} \rho^{p}_{\vec{\T}_b,\theta,\lambda;t}
\ = \
- \rho^{p}_{\vec{\T}_b,\theta,\lambda;t}
\
+ \lambda \,
\, \P\left(\eta^{p}_{\vec{\T}_b,\theta,\lambda;t}(v) = 0,
\ \# \left\{
u \in \Cal{N}_{\vec{\T}_b,v} : \eta^{p}_{\vec{\T}_b,\theta,\lambda;t}(u) = 1
\right\} \geq \theta\right),
\Eq(generatorvecTb)
)
$$
and noticing that $\vec{\T}_b$ has the special property that 
the random variables
$\eta^{p}_{\vec{\T}_b,\theta,\lambda;t}(u)$, 
$u \in \Cal{N}_{\vec{\T}_b,v}$,
are independent and have the same distribution as
$\eta^{p}_{\vec{\T}_b,\theta,\lambda;t}(v)$.
If the random variable $\eta^{p}_{\vec{\T}_b,\theta,\lambda;t}(v)$
were also independent of those, \equ(generatorvecTb) would reduce to 
the mean-field equation \equ(MF1), but this independence does not hold. 
In each one of the two proofs below we deal with this 
lack of independence in a different way.

\demo{Proof of Theorem 1 and its analogue for $\vec{\T}_b$}
We will prove that for arbitrary $b \geq 2$ and $\theta \geq 2$, \
if $\lambda_{\text{c}}^{\text{MF}}(b,\theta) < 1$, and
$\lambda > \lambda_{\text{c}}^{\text{MF}}(b,\theta)/
(1 - \lambda_{\text{c}}^{\text{MF}}(b,\theta))$, then
$$
p_{\text{c}}(\vec{\T}_b,\theta,\lambda)
\ \leq \ p_{\text{c}}^{\text{MF}}(b,\theta,\lambda/(\lambda + 1)) \ < \ 1.
\Eq(goalproofthm1)
$$
This suffices, since it obviously implies
$$
\lambda_{\text{c}} (\vec{\T}_b,\theta)
\ \leq \ \frac{\lambda_{\text{c}}^{\text{MF}}(b,\theta)}{
1 - \lambda_{\text{c}}^{\text{MF}}(b,\theta)} 
\ < \ \infty,
\Eq(goalproofthm1')
$$
and \equ(Thm11) and \equ(Thm12) follow respectively from 
\equ(goalproofthm1') and \equ(goalproofthm1), since the threshold 
$\theta$ contact process on $\T_b$ stochastically dominates the 
threshold $\theta$ contact process on $\vec{\T}_b$.

Only the first inequality in \equ(goalproofthm1) needs to be proved. 
For this purpose consider the arbitrary site $v$ that appears in 
\equ(generatorvecTb) and define the following event $E$:
$\eta^{p}_{\vec{\T}_b,\theta,\lambda;0}(v) = 0$, and
either there is no $U$ nor $D$ mark at $v$ between times 0 and $t$,
or else, the last such mark is a $D$ mark. By using the time reversibility
of Poisson processes, a standard computation gives
$$
\P(E)
\geq \frac{1}{\lambda + 1} \ (1-p).
\Eq(PE)
$$
The event $E$ is clearly independent
of the random variables
$\eta^{p}_{\vec{\T}_b,\theta,\lambda;t}(u)$, 
$u \in \Cal{N}_{\vec{\T}_b,v}$.
Since also
$E \subset \{\eta^{p}_{\vec{\T}_b,\theta,\lambda;t}(v) = 0\}$,
it follows from \equ(generatorvecTb), \equ(PE) and the observation after
\equ(generatorvecTb) that
$$
\align
\frac{d}{dt} \rho^{p}_{\vec{\T}_b,\theta,\lambda;t}
& \ \geq \
- \rho^{p}_{\vec{\T}_b,\theta,\lambda;t}
\
+ \lambda \,
\, \P\left(E,
\  \# \left\{
u \in \Cal{N}_{\vec{\T}_b,v}: 
\eta^{p}_{\vec{\T}_b,\theta,\lambda;t}(u) = 1
\right\} \geq \theta\right)                              \\
& \ = \
- \rho^{p}_{\vec{\T}_b,\theta,\lambda;t}
\
+ \lambda \,
\, \P\left(E\right)
\ \P\left(  \# \left\{
u \in \Cal{N}_{\vec{\T}_b,v}
: \eta^{p}_{\vec{\T}_b,\theta,\lambda;t}(u) = 1
\right\} \geq \theta\right)                              \\
& \ \geq \
- \rho^{p}_{\vec{\T}_b,\theta,\lambda;t}
\
+
\ \frac{\lambda}{\lambda + 1} \ (1-p)
\ \text{Bin}\left(b,\rho^{p}_{\vec{\T}_b,\theta,\lambda;t},\theta\right).
\teq(generator>)
\endalign
$$
It is convenient to rewrite \equ(generator>) as
$$
\frac{d}{dt} \rho^{p}_{\vec{\T}_b,\theta,\lambda;t}
\ \geq \
L(b,\theta,\lambda,p;\rho^{p}_{\vec{\T}_b,\theta,\lambda;t}),
\Eq(generator>L)
$$
where
$$
L(b,\theta,\lambda,p;x) \ = \
x \, H(b,\theta, \lambda/(\lambda+1);x) \ + \
\frac{\lambda}{\lambda + 1} \ \text{Bin}(b,x,\theta) \ (x-p),
\qquad x \in [0,1].
$$
When $\lambda_{\text{c}}^{\text{MF}}(b,\theta) < 1$ and
$\lambda > \lambda_{\text{c}}^{\text{MF}}(b,\theta)/
(1 - \lambda_{\text{c}}^{\text{MF}}(b,\theta))$, then
$\lambda/(\lambda+1) > \lambda_{\text{c}}^{\text{MF}}(b,\theta)$.
So $p_{\text{c}}^{\text{MF}}(b,\theta,\lambda/(\lambda + 1)) < 1$,
and for $p > p_{\text{c}}^{\text{MF}}(b,\theta,\lambda/(\lambda + 1))$
arbitrarily close to  \newline
$p_{\text{c}}^{\text{MF}}(b,\theta,\lambda/(\lambda + 1))$ we have
$H(b,\theta, \lambda/(\lambda+1);p) > 0$.
%for some $p > p_{\text{c}}^{\text{MF}}(b,\theta,\lambda/(\lambda + 1))$.
It follows then that also \newline
$L(b,\theta,\lambda,p;p) > 0$. We claim that
%and $H(b,\theta, \lambda/(\lambda+1);x)$ is strictly increasing in
%$x$ at $x=p_{\text{c}}^{\text{MF}}(b,\theta,\lambda/(\lambda + 1))$.
%It follows then that there is
%$p > p_{\text{c}}^{\text{MF}}(b,\theta,\lambda/(\lambda + 1))$
%such that $L(b,\theta,\lambda,p;p) > 0$. We claim that
$$
\inf_{t \geq 0} \, \rho^{p}_{\vec{\T}_b,\theta,\lambda;t}
\ \geq p > 0.
\Eq(inft)
$$
Indeed, set
$$
t_p = \inf\{t \geq 0 : \rho^{p}_{\vec{\T}_b,\theta,\lambda;t} < p\}.
$$
If \equ(inft) were false, we would have $t_p < \infty$.
Then by the continuity of
$\rho^{p}_{\vec{\T}_b,\theta,\lambda;t}$ in $t$, we would have
$\rho^{p}_{\vec{\T}_b,\theta,\lambda;t_p} = p$, and
$d/dt \, \rho^{p}_{\vec{\T}_b,\theta,\lambda;t} \leq 0$, at $t = t_p$.
But \equ(generator>L) implies
$d/dt \, \rho^{p}_{\vec{\T}_b,\theta,\lambda;t} \geq
L(b,\theta,\lambda,p;p) > 0$, at $t = t_p$.
This contradiction proves \equ(inft), which implies
$$
p_{\text{c}}(\vec{\T}_b,\theta,\lambda) \leq p.
$$
Since $p$ can be taken arbitrarily close to
$p_{\text{c}}^{\text{MF}}(b,\theta,\lambda/(\lambda + 1))$,
the proof of \equ(goalproofthm1) is complete.
\cqd
\enddemo

\demo{Proof of Theorem 6}
Applying Harris' inequality to \equ(generatorvecTb), we obtain
$$
\frac{d}{dt} \rho^{p}_{\vec{\T}_b,\theta,\lambda;t}
\ \leq \
- \rho^{p}_{\vec{\T}_b,\theta,\lambda;t}
\
+ \lambda \,
\, \P\left(\eta^{p}_{\vec{\T}_b,\theta,\lambda;t}(v) = 0\right)
\ \P\left( \# \left\{
u \in \Cal{N}_v : \eta^{p}_{\vec{\T}_b,\theta,\lambda;t}(u) = 1
\right\} \geq \theta\right).
%\Eq(Harris)
$$
From the observation after \equ(generatorvecTb), now
$$
\frac{d}{dt} \rho^{p}_{\vec{\T}_b,\theta,\lambda;t}
\ \leq \
- \rho^{p}_{\vec{\T}_b,\theta,\lambda;t}
+
\lambda \left(1-\rho^{p}_{\vec{\T}_b,\theta,\lambda;t}\right)
\text{Bin}\left(b,\rho^{p}_{\vec{\T}_b,\theta,\lambda;t},\theta\right),
\Eq(HarrisMF)
$$
which for $p > 0$ is equivalent to
$$
\frac{d}{dt} \log \left(\rho^{p}_{\vec{\T}_b,\theta,\lambda;t}\right)
\ \leq \
H\left(b,\theta, \lambda;\rho^{p}_{\vec{\T}_b,\theta,\lambda;t}\right).
\Eq(HarrisMFlog)
$$

If $p = 0$, then \equ(Thm42a) holds. Suppose that $p > 0$
and \equ(Thm41a) fails. Then \equ(HarrisMFlog)
implies that $\rho^{p}_{\vec{\T}_b,\theta,\lambda;t} \to 0$ as $t \to \infty$. 
The proof that \equ(Thm42a) holds then can be completed as the proof of 
Proposition 2. This shows that for each $\lambda$ and $p$ either 
\equ(Thm41a) or \equ(Thm42a) holds.

The statement about the set $\Cal{D}(\vec{\T}_b,\theta)$
follows now from the fact that the negation of \equ(Thm41a) is 
a ``finite-time condition'':
$$
\rho^{p}_{\vec{\T}_b,\theta,\lambda;t}
<
p_{\text{c}}^{\text{MF}}(b,\theta,\lambda),
\text{ \ for some $t \geq 0$}.
\Eq(negThm41a)
$$
If  \equ(negThm41a) holds for some $(\lambda,p)$, then, by continuity,
it also holds close to this point, with the same $t$. 
\cqd
\enddemo

\

\subheading{3. The regime of small $\lambda$}
\numsec=3
\numfor=1
In this section we will prove Theorems 2 and 3.  
We will abbreviate the notation, omitting 
$\T_b$, $\theta$ and $\lambda$ for instance in:
$$
\eta^p_{\T_b,\theta,\lambda;t} = \eta^p_t,
\qquad
\rho^p_{\T_b,\theta,\lambda;t} = \rho^p_t,
\qquad
\Cal{N}_{\T_b,v} = \Cal{N}_{v}.
$$

We will compare the threshold contact process on $\T_b$ with the 
similar process in which the spin of one of the neighbors of the 
root is frozen in the state 1. 
Recall the definition of $\Cal{L}$ from Subsection 1.4, and the corresponding
terminology and notation.
In our modified process the flip rates are as in 
the threshold $\theta$ contact process on $\T_b$, except for the site 
$-1$, where the spin is kept frozen in the state 1. We will start this
modified process with each site other than site $-1$ taking independently 
the value +1 with probability $p$, and being in state 0, otherwise.  
The notation 
%$\eta^{*,p}_{{\T}_b,\theta,\lambda;t} = 
$\eta^{*,p}_{t}$ 
will denote this process. Define also
$$
\sigma^{l,p}_t  \ = \ \P(\eta^{*,p}_{t}(l) = 1),
$$
and abbreviate $\sigma^{p}_t = \sigma^{0,p}_t$.
%and $\sigma_t = \sigma^{1}_t$.

For comparison, we will also consider the trivial threshold 0 contact 
process on $\T_b$, i.e., the process in which the spin of each site flips
independently of anything else, with 0 flipping to 1 at rate $\lambda$,
and 1 flipping to 0 at rate 1. Let $\pi^p_t$ be the probability that in 
this process a given
site is in state 1 at time $t$, when at time 0 this probability is 
set to $\pi^p_0 = p$. 
It is elementary that for every $p \in [0,1]$, $\pi^p_t \leq 
\pi^1_t \searrow \lambda/(\lambda + 1)$ as $t \to \infty$. In particular,
there is $\tilde t (\lambda)$ such that 
$$
\pi^p_t \ \leq \ \lambda, \text{ for $t \geq \tilde t (\lambda)$}.
\Eq(pitildeit)
$$

By attractivity, for any $p \in [0,1]$, $t \geq 0$ and $0 \leq l_1 \leq l_2$, 
$$
\rho^p_t \leq \sigma_t^{l_2,p} \leq \sigma_t^{l_1,p} \leq \pi^p_t.
\Eq(attractivitychain)
$$

\proclaim{Lemma 1}
For arbitrary $b \geq 2$, $\theta \geq 2$, $p \in [0,1]$ and $t \geq 0$,
$$
\frac{d}{dt} \, \sigma^{l,p}_t
\ \leq \
- \sigma^{l,p}_t 
\ + \
\lambda b \, \sigma^{l+1,p}_t, \quad l \geq 0.
\Eq(L11)
$$
And
$$
\frac{d}{dt} \, \sigma^{l,p}_t
\ \leq \
- \sigma^{l,p}_t
\ + \
\lambda \{ 
\pi^p_t b \sigma^{p}_t  +  \text{Bin}(b,\sigma^{p}_t,\theta)
\}, \quad l \geq 1.
\Eq(L12)
$$
\endproclaim

\demo{Proof}
We will use the following Terminology. For each site $v$ of $\T_b$, 
the $b$ sites in $\Cal{N}_{\vec{\T}_b,v}$ will be called forward 
neighbors of $v$ and the single site in 
$\Cal{N}_{{\T}_b,v} \backslash \Cal{N}_{\vec{\T}_b,v}$
will be called the backward neighbor of $v$.

 From the definition of $({\eta}^{*,p}_{t})$ and $\sigma^{l,p}_t$, 
$$
\align
\frac{d}{dt} \sigma^{l,p}_t
& \ = \
- \sigma^{l,p}_t
\
+ \lambda \,
\, \P\left({\eta}^{*,p}_{t}(l) = 0,
\ \# \left\{
u \in \Cal{N}_l : {\eta}^{*,p}_{t}(u) = 1
\right\} \geq \theta\right)                                \\
& \ \leq \
- \sigma^{l,p}_t
\
+ \lambda \,
\, \P\left(
\ \# \left\{
u \in \Cal{N}_l : {\eta}^{*,p}_{t}(u) = 1
\right\} \geq \theta\right).   
\teq(generatorL1)
\endalign
$$
%Note that site $l$ has $b$ neighbors that are at distance $l+1$
%from the site $-1$; we will call them outside neighbors of $l$. The 
%other neighbor of $l$ is at distance $l-1$ from the site $-1$; call 
%it inside neighbor of $l$. 

Inequality \equ(L11) follows from \equ(generatorL1) and
the observation that, since $\theta \geq 2$,
for the site $l$ to have at least $\theta$ occupied neighbors, it must have 
at least one occupied forward neighbor.

To derive \equ(L12) from \equ(generatorL1), we compare the process 
${\eta}^{*,p}_{t}$ with a further modified process 
in which the spins at the sites $-1$ and $l$ are both frozen in the state
1, while the spins at other sites evolve as in the threshold $\theta$ 
contact process. In this modified process, the spins at the neighbors
of $l$ evolve independently of each other. 
Note that by attractivity, the distribution
of $({\eta}^{*,p}_{t}(u))_{u \in \Cal{N}_l}$
is therefore stochastically dominated by a product measure in which the 
forward  neighbors of $l$ have probability $\sigma^{p}_t$ of being occupied,
while the inside neighbor of $l$ has probability $\pi^{p}_t$ of being
occupied. The probability in the r.h.s. of \equ(generatorL1) is now 
estimated from above by the probability that either the
backward neighbor of $l$ and at least one of its forward neighbors are both 
occupied (recall $\theta \geq 2$), or else that at least $\theta$ of its 
forward neighbors are occupied. 
\cqd
\enddemo

\proclaim{Lemma 2} For arbitrary $\theta \geq 2$ and $A \in (0,\infty)$,
there are $b^*,\delta^*,t^* \in (0,\infty)$,
such that if $b \geq b^*$ and $\lambda \leq A/b$, then for every $p \in [0,1]$
the following dichotomy holds. Either
$$
\sigma^{p}_{t}
\ \geq \
\frac{\delta^*}{b}, \text{ \ for all $t \geq t^*$},
\Eq(L21)
$$
or, for some $C \in (0,\infty)$,
$$
\sigma^{p}_{t}
\ \leq \
C \, e^{-0.6 \, t}, \text{ \ for all $t \geq 0$}.
\Eq(L22)
$$
\endproclaim

\noindent
{\it Remark:}
The exponential rate 0.6 in \equ(L22), could be replaced with any rate
smaller than 1, with minor modifications in the proof and a larger value
for $b^*$. In our proof of Theorem 2, all that we will need about this 
rate is that it is larger than 1/2. 

\demo{Proof}
For later convenience, we take $b^* = 9A^3$. 
We will use Lemma 1, and for this purpose we need to estimate $\pi^p_t$.
Under the assumptions in the lemma that we are proving, 
we have $\lambda \leq A/b \leq A/b^* = 1/(9A^2) =: \hat \lambda$. 
Let $(\hat{\pi}^p_t)_{t \geq 0}$ be defined in the same way as
$(\pi^p_t)_{t \geq 0}$, but with $\hat \lambda$ replacing $\lambda$.
Clearly $\pi^p_t \leq \hat{\pi}^p_t$. 
So, by \equ(pitildeit), 
there is $t^* = \tilde t (\hat \lambda) $ which depends on $A$, 
but not on $b$ or $\lambda$ (once they satisfy the conditions in 
the lemma), such that $\pi^p_t \leq \hat\lambda = 1/(9A^2)$,
for every $p \in [0,1]$ and $t \geq t^*$.

We use now the two inequalities in Lemma 1. The first one with 
$l=0$ and the second one with $l=1$. We suppose that 
$t \geq t^*$, so that we can use the estimate above on $\pi^p_t$. 
Since also $\lambda \leq A/b$, these inequalities read then
$$
\align
\frac{d}{dt} \, \sigma^{p}_t
& \ \leq \
- \sigma^{p}_t
\ + \
A \, \sigma^{1,p}_t,   \\
\frac{d}{dt} \, \sigma^{1,p}_t
& \ \leq \
- \sigma^{1,p}_t
\ + \
\frac{1}{9A}
\, \sigma^{p}_t  +  \frac{A}{b} \, \text{Bin}(b,\sigma^{p}_t,\theta).
\endalign
$$
Multiply the first of these inequalities by $1/(3\sqrt{A})$ and 
the second one by $\sqrt{A}$, and add the resulting inequalities
to obtain
$$
\frac{d}{dt} \, x_t 
\ \leq \
-\frac{2}{3} \, x_t  \ + \ 
\frac{A^{3/2}}{b} \, \text{Bin}(b,\sigma^{p}_t,\theta),
\qquad t \geq t^*,
\Eq(ddxsigma)
$$
where
$$
x_t \ = \ \frac{1}{3\sqrt{A}} \, \sigma^{p}_t  \ + \
           \sqrt{A} \, \sigma^{1,p}_t.
$$
By \equ(attractivitychain), 
$\sigma^{1,p}_t \leq \sigma^{p}_t$. Therefore we
obtain the following comparison between $x_t$ and $\sigma^{p}_t$:
$$
\frac{1}{3\sqrt{A}} \, \sigma^{p}_t
\ \leq \
x_t 
\ \leq \
\left\{
\frac{1}{3\sqrt{A}} \, + \, \sqrt{A}
\right\} \sigma^{p}_t
\Eq(compxsigma)
$$
From \equ(ddxsigma), the fact that 
$\text{Bin}(b,\sigma^{p}_t,\theta) \leq b^2 (\sigma^{p}_t)^2 / 2$
(since $\theta \geq 2$) 
and the first inequality in \equ(compxsigma),
$$
\frac{d}{dt} \, x_t
\ \leq \
-\frac{2}{3} \, x_t  \ + \
\frac{3}{2} A^{5/2}\, b \, (x_t)^2
\ =: \ G(x_t),
\qquad t \geq t^*.
$$
Note that $G(x) < 0$ for $0 < x < 4/(9A^{5/2}b)$, and 
$\lim_{x \searrow 0} G(x)/x = -2/3$.
One can use these facts, as in the proof of Proposition 2, to conclude that 
if there is some $\tilde t \geq t^*$ such that 
$$
x_{\tilde t} \ < \ \frac{4}{9 \, A^{5/2} \, b},
\Eq(touseappendix)
$$
then there is some $C' < \infty$ such that
$$
x_t \leq C' \, e^{-0.6 t}, \text{ for all $t \geq 0$}.
\Eq(L22x)
$$
But from \equ(compxsigma), the condition \equ(touseappendix) is
implied by 
$$
\sigma^p_{\tilde t} \ < \ \frac{4}{3A^2\,(1+3A)\,b},
$$
and \equ(L22x) implies \equ(L22) with $C = 3 \sqrt{A} \, C'$.
This completes the proof of the lemma, with 
$\delta^* = 4/(3A^2(1+3A))$.
\cqd 
\enddemo

\proclaim{Lemma 3} For arbitrary $\theta \geq 2$ and $A \in (0,\infty)$,
there is $b_1 \in (0,\infty)$,
such that if $b \geq b_1$ and $\lambda \leq A/b$, 
then for every $p \in [0,1]$ and $l \geq 0$,
$$
\rho_t^p \ \geq \ \sigma_t^{l,p} - \left(\frac{1}{\sqrt b}\right)^{l+1}
\Eq(sigmarho)
$$
\endproclaim

\demo{Proof}
To derive 
\equ(sigmarho)
we consider the discrepancies between the processes
$(\eta^{p}_{t})$ and $(\eta^{*,p}_t)$. We construct these two processes
using the same structure of Poisson marks, and 
for every site $v$ of $\T_b$ we set
$$
\delta_t(v)  \ = \
\eta^{*,p}_t (v) \ - \ \eta^{p}_{t} (v).
$$
In this fashion, $\delta_t$ takes the value 1 at the sites where the 
processes disagree at time $t$, and 0 at the other sites.

Observe that a
$D$ mark eliminates a discrepancy, while a $U$ mark will 
possibly create a discrepancy at a site $v$ only if at least one 
neighbor of $v$ has a discrepancy at the time of this mark.
Therefore $\delta_t \leq \zeta^{*,0}_t$, where $(\zeta^{*,0}_t)$ is 
a process in which the spin at the site $-1$ is frozen in the state 1
while all other spins are initially set to 0, and evolve with the 
following rules at each site $v \not = -1$:
\roster
\item"$\bullet$" 1 flips to 0 at Poisson $D$ marks at $v$.
\item"$\bullet$" 0 flips to 1 at Poisson $U$ marks at $v$
if and only if at least one of the spins in $\Cal{N}_v$ is in state 1 at 
the time of that mark.
\endroster
Note that these flip rules are those of the threshold 1 contact process
on $\T_b$, with infection parameter $\lambda$. 

The threshold 1 contact process is stochastically 
dominated by the contact process, in which
spins flip at rates:
\roster
\item"$\bullet$" 1 flips to 0 at rate 1.
\item"$\bullet$" 0 flips to 1 at rate $\lambda$ times the number 
of sites of $\Cal{N}_v$ in state 1 at time $t$.
\endroster
We denote by $(\xi^{*,0}_t)$ the process on $\T_b$, in which 
the spin at the site $-1$ is frozen in the state 1
while all other spins are initially set to 0, and then allowed to evolve 
according to the flip rates of this contact process. 

The chain of comparisons presented above implies that 
$$
\sigma_t^{l,p} - \rho_t^p 
\ = \ 
\P (\delta_t(l) = 1) 
\ \leq \
\P (\zeta^{*,0}_t(l) = 1)
\ \leq \
\P (\xi^{*,0}_t(l) = 1).
$$

Let $(\xi^{\{v\}}_t)$ denote the contact process started from the 
configuration in which only the site $v$ is occupied (and no spin
is frozen). Self-duality for the contact process implies
$$
\align
\P (\xi^{*,0}_t(l) = 1)
& \ = \
\P (\xi^{\{l\}}_s(-1) = 1, \text{ for some $s \in [0,t]$})  \\
& \ \leq \ 
\P (\xi^{\{l\}}_s(-1) = 1, \text{ for some $s \geq 0$})     \\
& \ = \
\P (\xi^{\{0\}}_s(l+1) = 1, \text{ for some $s \in [0,t]$}) 
\ =: \ u(l+1).
\endalign
$$
The function $u(\cdot)$ has played an important role in the study of the
contact process on $\T_b$. The proof of \equ(sigmarho) will be complete once 
we argue that under our hypothesis, 
$$
u(l) \ \leq \  \left(\frac{1}{\sqrt b}\right)^{l}.
$$
For this purpose we refer to results in Chapter 4 in Part I of 
[Lig2], where references
to the original contributions can be found. The contact process on $\T_b$, 
$b \geq 2$,
has two critical points $0 < \lambda_1(b) < \lambda_2(b) < \infty$. 
Theorem 4.1 in Part I of [Lig2] tells us that 
$\lambda_2(b) \geq 1/(2 \sqrt b)$. 
Therefore we can find $b_1$ so that 
$\lambda \leq A/b < \lambda_2(b)$, when $b \geq b_1$.
Display (4.49) of Part I of 
[Lig2] tells us that $u(l) \leq (\beta(\lambda))^l$,
where $\beta(\lambda) := \lim_{l \to \infty} (u(l))^{1/l}$. 
Finally Theorem 4.65 in Part I of [Lig2] tells us that 
%$\beta(\cdot)$ is non-decreasing and
$\beta(\lambda) \leq 1/\sqrt b$ when $\lambda \leq \lambda_2(b)$.
This completes the proof of \equ(sigmarho).
\cqd
\enddemo

\proclaim{Lemma 4} For arbitrary $\theta \geq 2$ and $A \in (0,\infty)$,
there are $b_0, \delta, \delta^* ,t^* \in (0,\infty)$,
such that if $b \geq b_0$ and $\lambda \leq A/b$, then for every $p \in [0,1]$
the following dichotomy holds. Either
\equ(L21) and \equ(Thm41) both hold, 
or else \equ(L22) and \equ(Thm42) both hold.
\endproclaim

\demo{Proof}
Let $b^*$, $\delta^*$ and $t^*$ be as in Lemma 2. 
We will take $b_0 \geq b^*$, so that under the hypothesis of the lemma 
that we are proving we know from Lemma 2 that 
either \equ(L21) or \equ(L22) holds.

Suppose first that \equ(L21) holds. Define
$$
\bar\rho^p \ = \ \liminf_{T \to \infty} \frac{1}{T} \int_0^T \rho^p_t dt,
\qquad
\bar\sigma^{l,p} \ = \ \liminf_{T \to \infty} \frac{1}{T} \int_0^T 
\sigma^{l,p}_t dt,
\qquad
\bar\sigma^{p} \ = \ \bar\sigma^{0,p}.
$$
 From \equ(L11) and $\lambda \leq A/b$, for each $l \geq 0$,
$$
A \, \frac{1}{T} \int_0^T \sigma^{l+1,p}_t dt 
\ \geq \
\frac{1}{T} \int_0^T \sigma^{l,p}_t dt 
\ + \ \frac{\sigma^{l,p}_T - \sigma^{l,p}_0}{T}.
$$
Hence, 
$$
A \, \bar\sigma^{l+1,p} \ \geq \ \bar\sigma^{l,p}.
$$
By induction in $l$ and \equ(L21), we obtain now
$$
\bar\sigma^{l,p} 
\ \geq \
\frac{\bar\sigma^{p}}{A^l}
\ \geq \
\frac{\delta^*}{A^l \, b}.
\Eq(sigmabarb)
$$
We use now Lemma 3, and for this suppose that $b \geq \max \{b^*, b_1\}$.
Then, from \equ(sigmarho) and \equ(sigmabarb), we obtain
$$
\bar\rho^p
\ \geq \
\bar\sigma^{l,p} - \left(\frac{1}{\sqrt b}\right)^{l+1}
\ \geq \
\frac{\delta^*}{A^l \, b} - \left(\frac{1}{\sqrt b}\right)^{l+1}
\ \geq \
\frac{\delta^*}{A^l \, b} 
\left( 
1 - \frac{A}{\delta^*} \left( \frac{A}{\sqrt b} \right)^{l-1}
\right).
$$
Taking $b_0 \geq \max \{b^*, b_1\}$ large enough, we have 
$A/\sqrt b \leq 1/2$, when $b \geq b_0$. Hence there is $\hat l$
such that 
$$
\bar \rho^p
\ \geq \
\frac{\delta^*}{2 A^{\hat l} b},
$$
for all $b \geq b_0$. We conclude that \equ(Thm41) holds then 
with $\delta = \delta^*/(2 A^{\hat l})$.

Suppose now that \equ(L22) holds. 
From the definition of $({\eta}^{p}_{t})$ and $\rho^{p}_t$,
$$
\align
\frac{d}{dt} \rho_t^p
& \ = \
- \rho_t^p
\
+ \lambda \,
\, \P\left({\eta}^{p}_{t}(0) = 0,
\ \# \left\{
u \in \Cal{N}_l : {\eta}^{p}_{t}(0) = 1
\right\} \geq \theta\right)                                \\
& \ \leq \
- \rho_t^p
\
+ \lambda \,
\, \P\left(
\ \# \left\{
u \in \Cal{N}_l : \eta^{p}_{t}(0) = 1
\right\} \geq \theta\right).
\teq(generatorTb)
\endalign
$$
We compare the process $(\eta^{p}_{t})$ with 
the modified process in which the spin at the root is frozen in the state
1, while the spins at other sites evolve as in the threshold $\theta$ 
contact process. In this modified process, the spins at the neighbors
of the root evolve independently of each other. 
Note that by attractivity, the distribution
of $({\eta}^{p}_{t}(u))_{u \in \Cal{N}_0}$
is therefore stochastically dominated by a product measure 
with density $\sigma^{p}_t$. Therefore,
$$
\frac{d}{dt} \rho^p_t
+ \rho^p_t    
\ \leq \
\lambda \, \text{Bin}(b+1,\sigma^{p}_t,\theta)
\ \leq \
\lambda \, \frac{(b+1)^2}{2} \, (\sigma^{p}_t)^2
\ \leq \
C' \, e^{-1.2 t},
$$
for some $C' < \infty$, where in the last step we used \equ(L22).
Multiplying both sides of this differential inequality by $e^t$ 
and integrating yields \equ(Thm42) (see the end of the proof of 
Proposition 2 for an identical estimate). 
\cqd
\enddemo

\demo{Proof of Theorem 2}
Lemma 4 established the claimed dichotomy. It also implied that 
under the hypothesis of Theorem 2,
%when $b \geq b_0$ and $\lambda \leq A/b$, 
\equ(Thm42) is equivalent to the negation of \equ(L21), i.e.,
$$
\sigma^{p}_{t}
\ < \
\frac{\delta^*}{b}, \text{ \ for some $t \geq t^*$}.
\Eq(L21c)
$$
The statement about the set $\Cal{D}_A(\T_b,\theta)$
follows then from the fact that 
if \equ(L21c) holds for some $(\lambda,p)$, then it also holds close
to this point (with the same $t$). (As in the proof of Theorem 6, this 
is a typical ``finite-time condition'' argument.)
\cqd
\enddemo

\demo{Proof of Theorem 3}
Since \equ(Thm31<) has already been proved, we only have to prove that
$$
\liminf_{b \to \infty} \, b \, \lambda_{\text{exp}} 
\ \geq \
\Phi_{\theta}.
\Eq(Thm31lb)
$$
For this purpose, let $A < \Phi_{\theta}$ and $\lambda = A/b$. We will 
show that then
$$
\rho^1_{\infty} \ \leq \ \frac{1}{b^{3/2}}, \quad \text{when $b$ is large.}
\Eq(violateThm41)
$$
From Theorem 2 and the remarks after that theorem, we know that this implies 
that, when $b$ is large, alternative \equ(Thm42) must hold and hence 
$\lambda \leq \lambda_{\text{exp}}$. 
Therefore $b \, \lambda_{\text{exp}} \geq A$, and since $A$ can be taken 
arbitrarily close to $\Phi_{\theta}$, \equ(violateThm41)
implies \equ(Thm31lb).

From the proof of Proposition 1, in the introduction, we know that if 
$p \geq \lambda / (\lambda + 1)$,
then $\rho^1_{\infty} \leq \rho^p_t$, for every $t \geq 0$. Therefore 
\equ(violateThm41) will follow once we show that 
$$
\inf_{t \geq 0} \,\rho^{\lambda}_{t} 
\ \leq \ \frac{1}{b^{3/2}}, \quad \text{when $b$ is large.}
\Eq(violateThm41p)
$$

To prove this claim, we 
use again \equ(generatorTb), but this time we
compare the process $(\eta^{p}_{t})_{t \geq 0}$ 
with the modified process
in which the spin at the root is frozen in the state
0, while the spins at other sites evolve as in the threshold $\theta$
contact process. We denote this modified process by 
%$(\eta^{\cdot,p}_t)_{t \geq 0}$.
$(\eta^{\odot,p}_t)_{t \geq 0}$.
Let $F_T$ be the event that the origin is vacant at time 0 and that between
time 0 and time $T$ there is no $U$ mark at the origin. 
Then, for $0 \leq t \leq T$,
$$
\align
\P\left(
\ \# \left\{
u \in \Cal{N}_l : \eta^{p}_t(u) = 1
\right\} \geq \theta\right)
& \ \leq \
\P\left(
\ \# \left\{
u \in \Cal{N}_l : \eta^{p}_t(u) = 1
\right\} \geq \theta , \ F_T \right) \ + \ \P((F_T)^c)  \\
& \ = \
\P\left(
\ \# \left\{
u \in \Cal{N}_l : \eta^{\odot,p}_t(u) = 1
\right\} \geq \theta , \ F_T \right) \ + \ \P((F_T)^c)  \\
& \ \leq \
\P\left(
\ \# \left\{
u \in \Cal{N}_l : \eta^{\odot,p}_t(u) = 1
\right\} \geq \theta \right) \ + \ 
p \ + \ \lambda T.
\endalign
$$
In the modified process $(\eta^{\odot,p}_t)_{t \geq 0}$,
the spins at the neighbors
of the root evolve independently of each other.
Note that by attractivity, the distribution
of $(\eta^{\odot,p}_t(u))_{u \in \Cal{N}_0}$
is therefore stochastically dominated by a product measure
with density $\rho^{p}_t$. Therefore, \equ(generatorTb) yields
$$
\frac{d}{dt} \rho^p_t
\ \leq \
- \rho^p_t  \ + \
\lambda \, \left( \text{Bin}(b+1,\rho^{p}_t,\theta) \ + \
p \ + \ \lambda T \right),
\qquad 0 \leq t \leq T.
$$
We will use this inequality with $p = \lambda = A/b$ and $T = b^{1/4}$.
We also change variables to $x_t = (b+1) \, \rho^{p}_t$.
The inequality above then implies
$$
\frac{d}{dt} x_t 
\ \leq \
- x_t  \ + \
A \, \text{Bin}\left(b+1,\frac{x_t}{b+1},\theta\right) \ + \
\frac{A}{b} \ + \
\frac{2 A^2}{b^{3/4}}, 
%A^2 / b  + A^2 / b^{3/4} \right),
\qquad 0 \leq t \leq b^{1/4}.
$$
If \equ(violateThm41p) were false, then 
there would be arbitrarily large $b$ for which
$$
x_t > \frac{1}{b^{1/2}}, \qquad 
\text{for all $t \geq 0$}.
%0 \leq t \leq b^{1/4}.
\Eq(notviolateThm41p)
$$ 
We would then have, for the $b$ for which \equ(notviolateThm41p) holds,  
$$
\frac{d}{dt} \log(x_t)
\ \leq \
- 1  \ + \
\frac{A}{x_t} \, \text{Bin}\left(b+1,\frac{x_t}{b+1},\theta\right) \ + \
\frac{A}{b^{1/2}}  \ + \
\frac{2 A^2}{b^{1/4}},
\qquad 0 \leq t \leq b^{1/4}.
\Eq(ddtlogxt)
$$
Thanks to \equ(bintopoisson), as $b \to \infty$,
the right hand side of this inequality, as a function of $x_t$, converges
uniformly to 
$$
F(x_t) = -1 \ + A \, \ \frac{\text{Poisson}(x_t,\theta)}{x_t}.
$$
Since $A < \Phi_{\theta}$, we have $\sup_{x > 0} F(x) = -C$, for some $C > 0$. 
Therefore \equ(ddtlogxt) yields, when $b$ is large, 
$$
\frac{d}{dt} \log(x_t)
\ \leq \
- \frac{C}{2}, \qquad 0 \leq t \leq b^{1/4}.
$$
This implies
$$
x_{b^{1/4}} \ \leq \ x_0 \, e^{-(C/2)b^{1/4}}
\ = \ 
A \, \frac{b+1}{b} \, e^{-(C/2)b^{1/4}},
%A \, \frac{b+1}{b} \, e^{-\frac{C}{2}b^{1/4}},
$$
which for large $b$ contradicts \equ(notviolateThm41p). 
This contradiction  proves \equ(violateThm41p), and completes the proof of 
\equ(Thm31lb).
\cqd
\enddemo

\

\subheading{4. The regime of large $\lambda$}
\numsec=4
\numfor=1
In this section we will prove Theorems 4 and 5.
Our lower bounds on $p_{\text{exp}}(\T_b,\theta, \infty)$ 
will be obtained by comparison with a bootstrap percolation model 
that we describe next. 

The continuous time bootstrap percolation model on the  
graph or oriented graph $G = (V,E)$
%${G} = (V, E)$, 
with threshold $\theta$ and infection parameter $\lambda$
can be defined by taking the threshold contact process 
$(\eta_{{G},\theta,\lambda;t})_{t \geq 0}$
on ${G}$, with same 
threshold $\theta$, and suppressing all the flips from 1 to 0. In 
other words, the bootstrap percolation process has flip rates
at $v \in V$ at time $t \geq 0$ given by
\roster
\item"$\bullet$" 1 flips to 0 at rate 0.
\item"$\bullet$" 0 flips to 1 at rate $\lambda$ in case there are at least
$\theta$ sites of $\Cal{N}_{{G},v}$ in state 1 at time $t$, and at
rate 0 otherwise.
\endroster
We will denote by $(\zeta^{\mu}_{{G},\theta,\lambda;t})_{t \geq 0}$
%@@ zeta = \bar{\eta} ?
the resulting process, 
started from a random distribution picked according
to law $\mu$ at time 0.
If one uses the same Poisson system of $D$ and $U$ marks to construct 
$(\eta^{\mu}_{{G},\theta,\lambda;t})$ and 
$(\zeta^{\mu}_{{G},\theta,\lambda;t})$, then, clearly 
$$
\eta^{\mu}_{{G},\theta,\lambda;t}
\ \leq \
\zeta^{\mu}_{{G},\theta,\lambda;t},
\quad \text{for all $t \geq 0$}.
\Eq(tcp<bp)
$$

It is also clear that $\zeta^{\mu}_{{G},\theta,\lambda;t}$ is increasing 
in time, and therefore has a limit, 
$\zeta^{\mu}_{{G},\theta,\lambda;\infty}$. 
It is also clear that $\zeta^{\mu}_{{G},\theta,\lambda;\infty}$ does 
not depend on $\lambda > 0$, and that it can be obtained by the 
following iteration. 
Let $S_0$ be the set of sites which at time 0 are in state 1. 
Recursively define then
$$
S_{n} \ = \ S_{n-1} \, \cup \, \{v \in (S_{n-1})^c :
\# \{\Cal{N}_{{G},v} \cap S_{n-1} \}  \geq \theta  \},
\quad n \geq 1.
$$
The sets $S_n$ increase, and their limit is
$\cup_n S_n = S_{\infty} = \zeta^{\mu}_{G,\theta,\lambda;\infty}$.
(This iteration is often taken as the definition of bootstrap
percolation in discrete time.)

 From \equ(tcp<bp), it follows that, for any $\lambda > 0$,
the sites that are vacant in 
%$\zeta^{\mu}_{{G},\theta,\lambda;\infty}$ are vacant in the process 
$S_{\infty}$ are vacant in the process 
$(\eta^{\mu}_{{G},\theta,\lambda;t})$ at all times. 

When $G = \T_b$, the observation in the last paragraph can be used 
as follows.
Consider the clusters of occupied sites in 
$S_{\infty}$, i.e., 
the connected components of the subgraph of $\T_b$ induced by the 
sites in $S_{\infty}$.
It is easy to see that since $\theta \geq 2$,
the sites that belong to finite clusters of
$S_{\infty}$
will eventually be in state 0 in the process
$\eta^{\mu}_{{G},\theta,\lambda;t}$. 
Suppose that $\mu$ is product measure with density $p$. 
If $S_{\infty}$ contains almost surely only finite clusters, then 
for any $\lambda > 0$ the process 
$(\eta^{p}_{{\T_b},\theta,\lambda;t})_{t \geq 0}$ dies out.
Therefore $p \leq p_c(\T_b,\theta, \infty)$. 
The next lemma provides a stronger conclusion under a stronger assumption.
In its statement and its proof, we will use the following terminology 
and notation for $\T_b$. The distance between two sites is the length
of the path that connects them. We will use $B_n$ for the ball of radius 
$n$ and center at the origin. The outside neighbors of a site $v$ are 
the neighbors of $v$ that are farther apart from the origin than $v$
(each site $v \not = 0$ has $b$ outside neighbors, and the origin has
$b+1$ outside neighbors). 
Denote by $R_n$ the event that the site $0$ and some site separated
from it by distance $n$ are in the same cluster of $S_{\infty}$.
\proclaim{Lemma 5}
Suppose that for bootstrap percolation on $\T_b$, started from product 
measure with density $p$, $\P(R_n)$ decays exponentially with $n$.
%we have $\P(R_n) \leq C_1 e^{c_2 n}$ for some $C_1,C_2 \in (0,\infty)$. 
Then $p \leq p_{\text{exp}}(\T_b,\theta, \infty)$.
\endproclaim

\demo{Proof}
Let 
$(\eta^{B_n}_{\T_b,\theta,\lambda;t})$ be the
threshold contact process started with the ball of radius $n$
around the origin fully occupied and all other sites vacant.
From the observations above,
$$
\P\left(\eta^{p}_{\T_b,\theta,\lambda;t}(0) = 1, \, (R_n)^c\right) 
\ \leq \
\P\left(\eta^{B_n}_{\T_b,\theta,\lambda;t}(0) = 1\right).
\Eq(RncBn)
$$

The process $(\eta^{B_n}_{\T_b,\theta,\lambda;t})$
is stochastically dominated by 
the process started from the same configuration, in which a spin 0
never flips and a spin 1 flips to 0, at rate 1, iff all its outside 
neighbors are in state 0. 
For this process, let $\tau_v$, $v \in B_n$, be the random amount of 
time needed for the spin at $v$ to flip to 0 after the moment when 
it became allowed to flip. Clearly the $\tau_v$, $v \in B_n$, are 
i.i.d., with exponential distribution with mean 1. A simple induction 
argument, starting from the sites at distance $n$ from the root, and 
moving inwards, shows that the root will then flip to 0 at the random 
time
$$
\max_{\pi \in \Pi_n} \ \sum_{v \in \pi} \, \tau_v,
$$
where $\Pi_n$ is the set of $(b+1) \, b^{n-1}$ paths from 0 to the sites
that are at distance $n$ from it.

We obtain therefore, from \equ(RncBn), and the hypothesis of the lemma
$$
\align
\P\left(\eta^{p}_{\T_b,\theta,\lambda;t}(0) = 1 \right)
& \ \leq \  
\P(R_n)
\ + \ 
\P\left(\eta^{p}_{\T_b,\theta,\lambda;t}(0) = 1, \, (R_n)^c\right)  \\
& \ \leq \
%\P(R_n)
C_1 \, e^{-C_2 n}
\ + \
\P\left( 
\max_{\pi \in \Pi_n} \ \sum_{v \in \pi} \, \tau_v \ \geq \ t
\right)  \\
& \ \leq \
%\P(R_n)
C_1 \, e^{-C_2 n}
\ + \
(b+1) \, b^{n-1} \,
\P\left(
\sum_{v \in \tilde \pi} \, \tau_v \ \geq \ t
\right),
\endalign
$$
where $C_1, C_2$ are positive finite constants and
$\tilde \pi$ is an arbitrary element of $\Pi_n$.
Taking $n = \lfloor \epsilon t \rfloor$, for some $\epsilon > 0$
small enough, a standard large deviation estimates for Poisson random 
variables (see, e.g., (A.1) in the Appendix of [KS]) shows that 
$$
\P\left( 
\sum_{v \in \tilde \pi} \, \tau_v \ \geq \ t
\right)
\ \leq \
C_3 \, e^{-C_4 n},
$$
where $C_3 \in (0,\infty)$ and $C_4$ is large enough that 
$e^{-C_4} < b$.
The last two displayed inequalities imply then that 
$\P\left(\eta^{p}_{\T_b,\theta,\lambda;t}(0) = 1 \right)$ decays 
exponentially with $t$, completing the proof. 
\cqd
\enddemo

Bootstrap percolation on homogeneous trees has been studied in 
%@@ what is there?
[CLR] and [BPP]. 
Below we could build on some of their estimates. Nevertheless,
for the reader's benefit, and at little extra cost, we will present 
a self-contained approach to our problem 
of estimating $\P(R_n)$, in order to use Lemma 5.
%regarding to the clusters of $S_{\infty}$. 

To study the bootstrap percolation process on $\T_b$, it is 
convenient to study also, as a tool, the bootstrap percolation
processes on its subgraph $\T^+_b$, induced by the following subset $V_b^+$ 
of vertices. The set $V_b^+$ is the minimal set of vertices of $\T_b$
with the properties that $0 \in V_b^+$ and if $v \in V_b^+$ then 
$\Cal{N}_{\vec{\T}_b,v} \subset V_b^+$.
We will also consider bootstrap percolation on the oriented graph 
$\vec{\T}^+_b$, which has as set of vertices also $V_b^+$, and 
defined then by 
$\Cal{N}_{\vec{\T}^+_b,v} = \Cal{N}_{\vec{\T}_b,v}$.

First we observe that bootstrap percolation on 
$\T^+_b$ and on $\vec{\T}^+_b$ 
are strongly related to each other in the following way.
If we start them from a same set of occupied sites, $S_0$, then
for $n \geq 0$, 
either both will have the root in $S_{n}$ or neither one
will have it. To see this, given a set of sites $R$ of $\T^+_b$,
say that a site $v \in R$ is 
hidden from the root in $R$ if there is another site 
$u \in R$ which belongs to the path which connects the root to $v$. 
Observe that if in the iteration which defines $S_{n}$ 
for either one of the two processes that we are considering we 
eliminate all the sites that are hidden from the root
in $S_{n-1}$, we do not
change the truth or falsehood of the statement that the root 
belongs to $S_n$. But with this modification, the sets
$S_n$ are the same for both processes.
%@@ see below
%(We clarify that below the sets $S_n$ are as originally defined,
%without the modification used in the argument in this paragraph.) 

It is easy to write down 
a recursion for the probability $p^+_n$ that the root
belongs to $S_n$ in the bootstrap percolation process on 
$\T^+_b$ or $\vec{\T}^+_b$, started from product distribution with 
density $p$. In the last paragraph we argued that $p^+_n$ is the 
same for both processes. Now, for the process on $\vec{\T}^+_b$,
the root will belong to $S_{n}$ in case it 
belongs to $S_0$, or in case it does not belong to $S_0$, but 
at least $\theta$ of the sites in 
$\Cal{N}_{\vec{\T}^+_b,0}$ are in $S_{n-1}$.
This observation and some obvious facts about the geometry of 
$\vec{\T}^+_b$ yield:
$$
p_{n}^+ 
\ = \
p \ + \ (1-p) \, \text{Bin}(b,p_{n-1}^+,\theta),
\Eq(BPrecursion)
$$
with initial condition $p_0^+ = p$.

The right-hand-side of \equ(BPrecursion) is an increasing continuous
function of $p_n^+$. Therefore 
$$
p^+_n \ \nearrow \  p^+_{\infty} \ := \
\inf \{ x > 0 : x = p + (1-p) \text{Bin}(b,x,\theta) \}.
\Eq(BPrecursioninfty)
$$
The limit $p^+_{\infty}$ is the probability that the root belongs 
to $S_{\infty}$ in this bootstrap percolation process on 
$\T_b^+$ or $\vec{\T}_b^+$.
Since $\theta \geq 2$,
$p^+_{\infty} \ \leq \
\inf \{ x > 0 : x = p + (b/2)^2 x^2\}.$
It is easy to use this observation to conclude that 
$$
p^+_{\infty} \searrow 0 \quad \text{as $p \searrow 0$}.
\Eq(ppinfty)
$$

The following concept will be used in the proof of Theorem 4. 
Consider bootstrap percolation on $G = (V,E)$, and let $W \subset V$
and $w \in W$. We will say that ``$w$ is eventually $W$-internally
occupied'' 
%by bootstrap percolation 
if $w$ becomes eventually occupied in the bootstrap percolation 
process restricted to $W$. To make the definition precise, set
$S_0^W = S_0 \cap W$, 
$$
S_{n}^W \ = \ S_{n-1}^W \, \cup \, \{v \in 
W\, \backslash \, S_{n-1}^W  \ : \
\# \{\Cal{N}_{{G},v} \cap S_{n-1}^W \}  \geq \theta  \},
\quad n \geq 1.
$$
We now say that $w$ is eventually $W$-internally occupied in case
$w \in S_{\infty}^W := \cup_n S_n^W$. 
%Note that the event that $w$ is eventually $W$-internally occupied 
Note that this event 
depends only on the initial configuration of occupied sites in $W$. 

\demo{Proof of Theorem 4}
We will use Lemma 5, and for this purpose we need to estimate $\P(R_n)$.
Consider bootstrap percolation on $\T_b$ and 
let $\{0 \leftrightarrow n \}$ denote the event that 
the sites $0$ and $n$ belong to the same cluster of $S_{\infty}$.
Note that this is the same as the event that the sites 
$0, 1, ..., n$ are all eventually occupied in this bootstrap
percolation process.  

We will use the definition in the last paragraph before this
proof in the case $G = \T_b$, 
%$W = V \, \backslash \, \Cal{L}$. 
$W = \Cal{L}^c$.
Note that the subgraph of
$\T_b$ induced by $W$ (i.e., the subgraph of $\T_b$ obtained by 
removing $\Cal{L}$ from the set of vertices, along with the edges 
incident to these vertices) is an infinite collection of copies 
of $\T_b^+$. 
The roots of these copies of $\T_b^+$ are 
neighbors to the sites in $\Cal{L}$,  with each site in $\Cal{L}$
being neighbor to $b-2$ of these roots.
For $j \in \Cal{L}$, 
define $X_j$ as the number of neighbors of the site $j$
which are eventually $W$-internally occupied.
When bootstrap percolation on $\T_b$ is started from product measure 
with density $p$,
it follows from the remarks above that each $X_j$ has a
binomial distribution corresponding to $b-2$ attempts each
with probability $p^+_{\infty}$ of success, where
$p^+_{\infty}$ is given by \equ(BPrecursioninfty). Clearly
the $X_j$ are also mutually independent.

For each $j \in \Cal{L}$
define a grade as follows. 
If the site $j$ is in state 1 at time 0, give this site grade A.
If not, give this site the grade according to:
also grade A if $X_j \geq \theta$,
grade B if $X_j = \theta - 1$, 
grade C if $X_j = \theta - 2$, 
grade F if $X_j \leq \theta - 3$.
The probability of obtaining grades A, B or C are then, 
respectively:
$$
\align
p_A  & \ = \  p \  + \  (1-p) \, \text{Bin}(b-2,p^+_{\infty}, \theta),
\\
p_B & \ = \ (1-p) \, \text{Bin}(b-2,p^+_{\infty}, \theta - 1),
\\
p_C & \ = \ (1-p) \text{Bin}(b-2,p^+_{\infty}, \theta - 2).
\teq(pABC)
\endalign
$$

Observe that if 
$\{0 \leftrightarrow n \}$ occurs, 
then the following must happen: 
\roster
\item"(i)" No site in $0,1,...,n$ can have grade F.
\item"(ii)" If the sites $j$ and $k$, with $0 \leq j < k \leq n$ 
both have grade C, then there must exist a site $i$ with $j < i < k$ 
with grade A.
\endroster
For a given realization of the process, denote by 
$n_A$, $n_B$ and $n_C$, respectively, the number of sites
in $\{0, ..., n\}$ which receive grades $A$, $B$, and $C$.
Then (i) implies
$$
n_A + n_B + n_C = n+1,
\Eq(sumn)
$$
while (ii) implies
$$
n_A \geq n_C - 1.
\Eq(AC)
$$
Since there are $3^{n+1}$ ways to assign grades A, B and C to the
sites in $\{0, ..., n\}$, it follows that
$$
\align
\P(0 \leftrightarrow n) 
& \ \leq \
3^{n+1} \max_{
\Sb
n_A,n_B,n_C
\\
n_A + n_B + n_C = n+1
\\
n_A \geq n_C - 1
\endSb}
(p_A)^{n_A} (p_B)^{n_B} (p_C)^{n_C}  \\
& \ \leq \
3^{n+1} \max_{
\Sb
n_A,n_B,n_C
\\
2 n_A + n_B \geq n
\endSb}
(p_A)^{n_A} (p_B)^{n_B} (p_C)^{n_C}  \\
& \ \leq \
3^{n+1} \max_{
\Sb
n_A,n_B
\\
n_A + n_B \geq n/2
\endSb}
(p_A)^{n_A} (p_B)^{n_B} \\
& \ \leq \
3^{n+1} \ (\max\{p_A,p_B\})^{n/2}.
\teq(P0nsimple)
\endalign
$$

 From the geometry of $\T_b$ and \equ(P0nsimple) we obtain
$$
\P(R_n) \ \leq \ (b+1) b^{n-1} \P(0 \leftrightarrow n)
\ \leq \
3 \, \frac{b+1}{b} \, 
\left[ 
3b \, (\max\{p_A,p_B\})^{1/2}
\right]^n.
$$
But \equ(ppinfty) and \equ(pABC) imply that for $p>0$ small enough,
$\max\{p_A,p_B\} < 1/(3b)^2$, and hence
$$
\P(R_n) \to 0, \quad \text{exponentially fast as $n \to \infty$}.
$$
Lemma 5 now implies 
$p_{\text{exp}}(\T_b,\theta, \infty) \geq p > 0$.
\cqd 
\enddemo

We turn now to the proof of Theorem 5. The origin of the exponent 
$\theta / (\theta - 1)$ there is the following. 
$\text{Bin}(b, \gamma/b^{\alpha},\theta)$ is of order $1/b^{\alpha}$ for 
large $b$ iff $\alpha = \theta / (\theta - 1)$.
The precise version of this statement that we need below
is the following one, which can be checked by elementary computations.
For arbitrary $\bar \gamma > 0$,
$$
\frac{\theta !}{\gamma^{\theta}} \, 
b^{\theta/(\theta-1)} \,
\text{Bin}(b,\gamma/b^{\theta/(\theta-1)},\theta) \ \to \
1,
%\frac{\gamma^{\theta}}{\theta !},
\quad \text{as $b \to \infty$},
\Eq(limb(theta/)
$$
uniformly in $\gamma \in (0, \bar \gamma]$.
This has the following consequence for the mean-field model.
For arbitrary $\bar \gamma > 0$,
$$
H(b,\theta,\lambda;\gamma/b^{\theta/(\theta-1)}) \ \to \
-1 \ + \ \lambda \ \frac{\gamma^{\theta-1}}{\theta !},
\quad \text{as $b \to \infty$},
$$
uniformly in $\gamma \in (0,\bar \gamma]$.
Therefore, for arbitrary $\lambda > 0$,
$$
b^{\theta/(\theta-1)} \, 
p_{\text{c}}^{\text{MF}}(b,\theta,\lambda) \ \to \ 
(\theta ! / \lambda)^{1/(\theta - 1)},
\quad \text{as $b \to \infty$}.
\Eq(pcbinftyMF)
$$
\demo{Proof of Theorem 5}
From \equ(Thm12), in Theorem 1, 
$$
p_{\text{c}}(\T_b,\theta,\infty) \ \leq \ 
p^{\text{MF}}_{\text{c}}(b,\theta,1).
$$
Combined with \equ(pcbinftyMF), this implies 
$$
\limsup_{b \to \infty} \  b^{\theta/(\theta-1)} \,
p_{\text{c}}({\T}_b,\theta,\infty)
\ \leq \
(\theta !)^{1/(\theta - 1)},
$$
which provides the upper bound in \equ(Thm32).

The proof of the lower bound in \equ(Thm32), 
$$
\liminf_{b \to \infty} b^{\theta/(\theta-1)} p_c(\T_b,\theta,\infty)
\ > \ 0,
\Eq(lastgoal)
$$
builds on the proof of Theorem 4.
We will use the same notation as in that proof.
Define also $\gamma = p \, b^{\theta/(\theta - 1)}$,
$\gamma_n = p^+_n \, b^{\theta/(\theta - 1)}$, and
$\gamma_{\infty} = p^+_{\infty} \, b^{\theta/(\theta - 1)}$.
The recursion \equ(BPrecursion) implies
$$
\gamma_{n} \ \leq \ \gamma \ + \ b^{\theta/(\theta - 1)} \, 
\text{Bin}(b,\gamma_{n-1}/b^{\theta/(\theta-1)},\theta)
\Eq(BPrecurssiongamma)
$$
Set
$$
\bar \gamma \ = \ 
\inf \{ x > 0 : x = \gamma + 2 x^{\theta} / \theta!  \}. 
%\Eq(balphapinfty)
$$
(In this definition the factor 2 is arbitrary; any number larger than 1
could be used instead.) 
Since $\theta \geq 2$, for small $\gamma > 0$ we have 
$0 < \bar \gamma < \infty$, and
$$
\bar \gamma \ = \ \gamma + 2 \frac{(\bar \gamma)^{\theta}}{\theta!}.
\Eq(bargamma)
$$ 
Moreover, similarly to \equ(ppinfty),
$$
\bar \gamma \searrow 0 \quad \text{as $\gamma \searrow 0$}.
\Eq(gammagammainfty)
$$
Our next goal is to prove that when $b$ is large
$$
p_{\infty}^+ \ \leq \ \frac{\bar \gamma}{b^{\theta/(\theta-1)}}
\Eq(balphapinfty)
$$
 From \equ(limb(theta/), there is $\bar b$ such that for 
$0 \leq y \leq \bar \gamma$ and $b \geq \bar b$,
$$
b^{\theta/(\theta-1)} \, 
\text{Bin}(b,y/b^{\theta/(\theta-1)},\theta)
\ \leq \
\frac{2y^{\theta}}{\theta !}.
\Eq(unifgammabar)
$$
Since $\gamma + 2 x^{\theta} / \theta!$ is increasing in $x$, 
when $y < \bar \gamma$, we have 
$\gamma + 2 y^{\theta} / \theta! \ \leq \ \bar \gamma$. 
Note that $\gamma \leq \bar \gamma$. Hence \equ(BPrecurssiongamma) and 
\equ(unifgammabar) imply, by induction on $n$, that 
$\gamma_n \leq \bar \gamma$.
Therefore $\gamma_{\infty} \leq \bar \gamma$ and 
\equ(balphapinfty) follows.

Combining \equ(pABC), \equ(balphapinfty) and
\equ(gammagammainfty) we obtain, when $b$ is large, 
$$
\align
p_A & \ \leq \ \frac{\gamma_A}{b^{\theta/(\theta-1)}},  \\
p_B & \ \leq \ \frac{\gamma_B}{b},  \\
p_C & \ \leq \ \frac{\gamma_C}{b^{(\theta - 2)/(\theta-1)}},
\teq(pABCgammaABC)
\endalign
$$
where $\lim_{\gamma \searrow 0} \gamma_A = 
\lim_{\gamma \searrow 0} \gamma_B = 
\lim_{\gamma \searrow 0} \gamma_C = 0$.
Note that \equ(pABCgammaABC) implies the following technical estimate:

$$
\max\{(p_A)^2,(p_B)^2,p_Cp_A \} \leq \ \frac{\gamma'}{b^2},
\Eq(maxpairs)
$$
where $\lim_{\gamma \searrow 0} \gamma' = 0$.
This estimate is useful in combination with the following one:
$$
\P(0 \leftrightarrow n)
\ \leq \
3^{n+1} \, \left( 
\max\{(p_A)^2,(p_B)^2,p_Cp_A \}
\right)^{(n/2) - 1}.
\Eq(P0nelaborate)
$$
To prove this inequality, one can match pairs of sites in 
$\{0,..., n\}$ in the following way. Recall that if the event 
$\{0 \leftrightarrow n \}$ happens, then the facts (i) and 
(ii) in the proof of Theorem 4 must happen. From fact (ii) we know that each 
site which receives a grade C is followed eventually by a 
site with grade A, except possibly for the last site with 
grade C. Pair each site with a grade C with the first site 
with grade A after it, leaving possibly one unmatched site
with grade C. Considering the sites with grade A which are 
unmatched to any site with grade C, we match the first of these
sites to the the second one, the third to the fourth, etc, 
leaving at most one unmatched site with grade A. Finally we 
match the first site with grade B to the second such site, 
the third site with grade B to the fourth such site, etc, 
leaving at most one unmatched site with grade B. Since the
number of sites in $\{0,..., n\}$ is $n+1$ and there are at 
most 3 unmatched sites, the number of matched pairs is
at least $((n+1)-3)/2 = (n/2) - 1$. The estimate 
\equ(P0nelaborate) now follows from the fact that the number 
of ways to assign grades A, B and C to the sites in $\{0, ..., n\}$,
is $3^{n+1}$. 

Combining \equ(maxpairs) with \equ(P0nelaborate), we obtain, when $b$ is large, 
$$
\P(R_n) \ \leq \ (b+1) b^{n-1} \P(0 \leftrightarrow n)
\ \leq \
(b+1) b^{n-1} 3^{n+1} \, \left( \frac{\gamma'}{b^2} \right)^{(n/2) - 1}
\ = \
27 (b+1)b \, \left(3 \sqrt{\gamma'}\right)^{n-2}.
$$
By taking $\gamma > 0$ sufficiently small, we can make $3 \sqrt{\gamma'} < 1$.
Then
$$
\P(R_n) \to 0, \quad \text{exponentially fast as $n \to \infty$}.
$$
Lemma 5 now implies
$ b^{\theta/(\theta-1)} \, p_{\text{exp}}(\T_b,\theta, \infty) \geq 
b^{\theta/(\theta-1)} \, p  = \gamma > 0$.
This proves \equ(lastgoal).
\cqd
\enddemo

\bigskip
\centerline{{\bf References}}
\bigskip

\item {\bf [BPP]}
{Balogh, J., Peres, Y., Pete, G.}
Bootstrap percolation on infinite trees and non-amenable groups.
{\it Combinatorics, Probability and Computing} (to appear).
\medskip

\item {\bf [BG]}
{Bramson, M. and Gray, L.} (1992) 
A useful renormalization argument.
In {\it Random walks, Brownian motion and interacting
particle systems. Festschrift in honor of Frank Spitzer}.
Durrett, R. and Kesten, H., editors. Birkh\"auser. 
pp. 113 - 152.
\medskip

\item {\bf [CLR]} 
{Chalupa, J., Leath, P. L., Reich, G. R.} (1979) 
Bootstrap percolation on a Bethe lattice. 
{\it Journal of Physics C} {\bf 12}, L31 -- L35. 
\medskip

\item {\bf [Chen1]} 
{Chen, H.-N.} (1992) On the stability of a population growth model with 
sexual reproduction on $\Z^2$. 
{\it Annals of Probability} {\bf 20}, 232--285. 
\medskip

\item {\bf [Chen2]} 
{Chen, H.-N.} (1994) On the stability of a population growth model with 
sexual reproduction on $\Z^d$, $d \ge 2$. 
{\it Annals of Probability} {\bf 22}, 1195--1226. 
\medskip

\item {\bf [DG]}
{Durrett, R. and Gray, L.} (1990)
Some peculiar properties of a particle system with sexual reproduction.
{\it Unpublished manuscript}.
\medskip

\item {\bf [KS]} 
{Kesten, H. and Schonmann, R. H.,} (1995)  
On some growth models with a small parameter. 
{\it Probability Theory and Related Fields} {\bf 101}, 435--468. 
\medskip

\item {\bf [Lig1]} 
{Liggett, T. M.,} (1985)  
Interacting particle systems. 
{\it Springer}, New-York, Berlin.
\medskip

\item {\bf [Lig2]}
{Liggett, T. M.,} (1999)
Stochastic Interacting Systems: Contact, Voter and Exclusion Processes.
{\it Springer}, New-York, Berlin.
\medskip

\item {\bf [Toom]}
{Toom, A. } (1974)
Nonergodic multidimensional systems of automata.
{\it Problems of information transmission} {\bf 10}, 239-246.
\medskip

\bigskip

\vfill\eject

Luiz Renato Fontes

Instituto de Matem\' atica e Estat\' \i stica 
 
Universidade de S\~ ao Paulo 

Rua do Mat\~ ao, 1010 - Cidade Universit\' aria

05508-090 \,\,  S\~ ao Paulo, SP 

Brazil

lrenato\@ime.usp.br

\
\
\

Roberto H. Schonmann

Mathematics Department

University of California at Los Angeles

Los Angeles, CA 90024

U.S.A.

rhs\@math.ucla.edu

\bye

\bye